\documentclass[12pt]{article}
\usepackage{latexsym}
\usepackage{amsfonts}
\usepackage{amsmath}
\usepackage{amsthm}
\usepackage{amssymb}

\newtheorem{theorem}{Theorem}[section]
\newtheorem{lemma}[theorem]{Lemma}

\newtheorem{proposition}[theorem]{Proposition}

\newtheorem{corollary}[theorem]{Corollary}
\newtheorem{remark}[theorem]{Remark}

\newtheorem{assumption}{Assumption}
\newtheorem{example}[theorem]{Example}
\numberwithin{equation}{section}
\renewenvironment{proof}{\textbf{Proof:}}{\par}

\usepackage[dvips]{color}
\usepackage{mathtools}
\mathtoolsset{showonlyrefs}
\usepackage{mathrsfs}
\usepackage{comment}
\usepackage{hyperref}
\usepackage{xcolor}

\def\R{\mathbb{R}}
\def\P{\mathbb{P}}
\def\E{\mathbb{E}}
\def\T{\mathbb{T}}
\def\FF{\mathcal{F}}
\def\LL{\mathcal{L}}
\def\S{\mathcal{S}}
\def\1{\mathbf{1}}
\def\mP{{\mathbf P}}
\def\mE{{\mathbf E}}
\def\d{{\rm d}}
\def\cH{\mathcal{H}}
\def\cD{\mathcal{D}}
\def\cI{\mathcal{I}}
\def\X{\mathbb{X}}
\def\Y{\mathbb{Y}}

\def\cL{{\cal L} }

\def\N{\mathcal{N}}

\def\R{\mathbb{R}}

\title{Large deviations and almost sure convergence for
the extremes of branching L\'evy processes}
\author{
	{\bf Runjia Luo}\thanks{R. Luo:  School of Mathematical Sciences, Capital Normal
		University,  Beijing, 100048, P.R. China.
		Email: {\texttt 2230502114@cnu.edu.cn}}
	\quad
	{\bf Yan-Xia Ren}\thanks{Y.-X. Ren:  LMAM School of Mathematical Sciences \& Center for
		Statistical Science, Peking
		University,  Beijing, 100871, P.R. China.
		Email: {\texttt yxren@math.pku.edu.cn}}
	\quad
	{\bf Renming Song}\thanks{R. Song: Department of Mathematics,
		,
		Urbana, IL 61801, U.S.A.
		Email: {\texttt rsong@illinois.edu}}
	\quad
	{\bf Rui Zhang*}\thanks{R. Zhang*:  corresponding author. School of Mathematical Sciences \& Academy for Multidisciplinary Studies, Capital Normal
		University,  Beijing, 100048, P.R. China.
		Email: {\texttt zhangrui27@cnu.edu.cn}}}		
\date{}
\begin{document}
	\maketitle
	\begin{abstract}
		In this paper, we investigate the asymptotic behavior  of  supercritical branching Markov  processes $\{\mathbb{X}_t, t \ge0\}$
		whose spatial motions are L\'evy processes with regularly varying tails.
		Recently, Ren et al. [Appl. Probab. 61 (2024)] studied the weak convergence
		of the extremes  of $\{\mathbb{X}_t, t \ge0\}$.
		In this paper, we establish the large deviation
		of $\{\mathbb{X}_t, t \ge0\}$
		as well as some almost sure convergence
		results of the maximum of $\mathbb{X}_t$.
	\end{abstract}
	

	\section{Introdution}
	
	\subsection{Model and notation}

	A  branching L\'evy process on $\R$ is a continuous-time particle system which can be described as follows.
	The system begins at time
	$t=0$ with a single particle located at $x$ which moves according to a branching L\'evy process $\{\xi_t,\mP_x\}$ with L\'evy exponent $\psi(\theta)$.
	After an independent exponential time with parameter $\beta$, the initial
	particle dies and gives birth to $k$ new particles with probability  $p_k$, $k\ge0$.
	Each new particle moves according to the L\'evy process $\xi$ starting from the position of its parent's death, and branches independently with the same branching rate $\beta$ and offspring distribution $\{p_k, k\ge 0\}$.
All particles, once born, evolve independently of one another.
The expectation with respect to $\mP_x$ will be denoted by $\mE_x$. We write $\mP:=\mP_0$ and $\mE:=\mE_0$.
	
	We label each particle using the classical Ulam-Harris labeling system. We denote by $\mathbb{T}$ the set of all particles in the tree and use $o$  represent the root of the tree.  For any $u\in \mathbb{T}$, let $N^u$ be the number of the offspring of $u$ and $\tau_u$  denote the lifetime of $u.$ Then $\{\tau_u,u\in\mathbb{T}\}$ are
i.i.d. exponential random variables of parameter $\beta$.

 We will use the notation $u<v$ to mean that $u$ is an ancestor of $v$, and $u\le v$ that either $u<v$ or $u=v$. For any $u\in \mathbb{T}$, let $I_u^0:=\{v\in\mathbb{T}:v<u\}$ and $I_u:=\{v\in\mathbb{T}:v\le u\}$.	Let $b_u$ and $\sigma_u$ be the birth time and death time of $u$ respectively.  It is clear that
 $$b_u=\sum_{v\in I_u^0}\tau_v, \qquad  \sigma_u=b_u+\tau_u.$$
 Let	$\mathcal{L}_t$  be  the set of all  particles alive at time $t$. Then $u\in\mathcal{L}_t$ means $b_u\le t<\sigma_u$.
 For $u\in\mathbb{T} \mbox{ with } b_u\le t$, define $\tau_{u, t}=\sigma_u\wedge t-b_u$.
 For any $u\in\mathcal{L}_t$, let $\xi_t^u$ be the position of $u$ at time $t$.
 	The branching L\'evy process $\{\mathbb{X}_t:t\ge0\}$ is the measure-valued process  defined by
	$$ \mathbb{X}_t:=\sum_{u\in\mathcal{L}_t}\delta_{\xi^u_t}. $$
	We use $\P_x$ to denote the law of the branching L\'evy process when the initial particle starts at position $x$. The expectation with respect to $\P_x$ will be denoted by $\E_x$. We write $\P:=\P_0$ and $\E:=\E_0$.
		Let $\{\mathcal{F}_t\}$ be the natural filtration of $\mathbb{X}$ and
	$$
	\mathcal{F}_t^{\mathbb{T}}:=\sigma(\{N^u: u\in\mathbb{T} \mbox{ with } \sigma_u\le t\}\cup\{\tau_{u,t}: u\in\mathbb{T} \mbox{ with } b_u\le t\}).
	$$

	In this paper, we study supercritical branching L\'evy processes,
that is to say,  we always assume  that $\mu:=\sum_k kp_k>1.$
		Let $\mathcal{S}$ be survival event. Then $\P_x(\mathcal{S})>0$ does not depend on the location $x$ of the initial particle.
The extinction probability $q:=\P(\mathcal{S}^c)$
is the unique root in the interval $[0,1)$ of the equation
$f(s)=s$, where $f(s):=\sum_{k}p_k s^k$.
	For more details, see \cite[Section III. 4]{Athreya-Ney}.
		For any $x\in \R$, we define $\P^*_x(\cdot):=\P_x(\cdot|\mathcal{S})$  and denote the corresponding expectation by $\mathbb{E}^*_x$.
		We write $\P^*:=\P^*_0$ and $\E^*:=\E^*_0$.

Recently, many people studied the extreme of branching L\'evy processes defined by
	\begin{align} \label{R_t}
		R_t := \max_{u \in \LL_t}\xi^u_t, \quad t>0.
	\end{align}
Here we use the convention that $\max\emptyset =-\infty.$
Among these, branching Brownian motions garnered the most attention.
For branching Brownian motions,
	Bramson  \cite{Bramson78}   (see also \cite{Bramson}) proved
	that,  under some moment conditions on the offspring distribution,
	$\P(R_t-m(t)\le x)\to 1-w(x)$
	as $t\to\infty$ for all $x\in \R$, where $m(t)=\sqrt{2}t-\frac{3}{2\sqrt{2}}\log t$ and  $w(x)$ is a traveling wave solution.
	For the large deviation  of  $R_t$,
	Chauvin and Rouault \cite{ Chauvin88,Chauvin} studied the
	asymptotic behavior of
	$P(R_t>\sqrt{2}\delta t)$ for $\delta \ge 1$.
	Derrida and Shi \cite{DS1, DS}
	studied the lower large deviation of $R_t$, i.e, the asymptotic behavior of $\frac1{t}\log P(R_t\le \sqrt{2}\delta t)$ for $\delta <1$, and found that the rate function has a phase transition at $1-\sqrt{2}$.
	Subsequently, Chen, He and Mallein \cite{CHM} studied the
	asymptotic behavior of
	$P(R_t\le \sqrt{2}\delta t)$ for $\delta<1.$
	Recently, \cite{ABBS, ABK12,ABK} studied the extremal processes of branching Brownian motions.

	In this paper we study the case when  the spatial motion
	is a L\'evy process with regularly varying L\'evy exponent.
	We now state our assumptions.
	
	Let $Z_t$ be the total number of the particles alive at time $t$. It is well known that $\{Z_t:t\ge 0\}$ is a continuous time Galton-Watson process.
		The following two quantities will play important roles in this paper:
	\begin{align}\label{d:lambda&vartheta}
	\lambda:=\beta(\mu-1), \quad \vartheta:=\int_0^\infty  e^{-\lambda r}\P(Z_r>0)\d r.
	\end{align}
	By the Markov property and the branching property,  the process $\{e^{-\lambda t}Z_t, t\ge0\}$ is a non-negative martingale with respect to $ \{\mathcal{F}^{\mathbb{T}}_t\}$. Thus it has an almost sure limit
	$$\lim_{t\to\infty}e^{-\lambda t}Z_t
	=:W.$$
	It is well known
	that $W$ is non-degenerate if and only if the following
	$L\log L$ criterien holds:
	\begin{assumption}\label{assum1}
		$\sum\limits_{k=1}^{\infty}(k\log k p_k) < \infty. $
	\end{assumption}
	\noindent Moreover,
under
{\bf Assumption}
\ref{assum1},
$\P(W>0)=\P(\mathcal{S})>0.$ For more details,
	see \cite[Section III.7]{Athreya-Ney}.

	We will always assume that the spatial motion satisfies:

	\begin{assumption}\label{assum2}
		There exist $\alpha \in(0,2)$, a complex constant $c_*$ with ${\rm Re}(c_*)>0$ and a function $L(x) : \R_+ \to \R_+$ slowly varying at $\infty$ such that $\psi(\theta) \sim -c_*\theta^{\alpha}L(\theta^{-1}) $  as $\theta \to 0_+$.
	\end{assumption}

Strictly $\alpha$-stable processes satisfy {\bf Assumption} \ref{assum2}.
By using the tables of complete Bernstein functions in \cite{SSV}, we can come up a lot of subordinate Brownian motions satisfying {\bf Assumption} \ref{assum2}.
Further discussions of
{\bf Assumption}
\ref{assum2} can be found in the Appendix.

In the Appendix, we will show that, under
{\bf Assumption}
\ref{assum2}, the function
$e^{-c_*|\theta|^\alpha}, \theta\in \R$,
is the characteristic function of an $\alpha$-stable random variable with L\'evy measure $v_\alpha$, where
	\begin{equation}\label{def-v-alpha}
		v_\alpha(dx)=q_1 x^{-1-\alpha}{\bf 1}_{(0,\infty)}(x)\d x+q_2 |x|^{-1-\alpha}{\bf 1}_{(-\infty,0)}(x)\d x,
	\end{equation}
	with $q_1$ and $q_2$ being nonnegative numbers,  uniquely determined by the following equation:
	if $\alpha\neq 1$
	$$c_*=\alpha\Gamma(1-\alpha)\left(q_1e^{-i\pi\alpha/2}+q_2e^{i\pi\alpha/2}\right),$$
	and if $\alpha=1$
	$$q_1=q_2=\Re(c_*)/\pi.$$
It has been proved in  \cite[Remark  2.1]{YRR} that for any $s>0$,
\begin{equation}\label{24}
		\mP (\xi_s\ge x)\sim  \frac{q_1}{\alpha} s x^{-\alpha}L(x),\quad \mP (\xi_s\le -x)\sim  \frac{q_2}{\alpha} s x^{-\alpha}L(x),\quad x\to\infty,
\end{equation}
	that is,
	$\xi_s$
	has regularly varying tails. To ensure the right tail of $\xi_s$ is regularly varying, we always assume that
	\begin{assumption}\label{assum3}
		$q_1>0$.
	\end{assumption}

In this paper, we always  assume that
{\bf Assumptions}
\ref{assum1}-\ref{assum3}
hold, and that $\alpha\in (0, 2)$, $c_*$ and $L$ are as specified in
{\bf Assumption}
 \ref{assum2}.

	The following variant of $\vartheta$ will also play a role later in this paper:
	\begin{align}\label{d:varthetastar}
	\vartheta^*:=
	\dfrac{q_1}{\alpha}\int_0^\infty  e^{-\lambda r}\P(Z_r>0)\d r
=\dfrac{q_1}{\alpha} \vartheta.
	\end{align}
	
	\bigskip
	
	Put  $\overline{\R}_0=[-\infty,\infty]\setminus\{0\}$
	with the topology generated by the set $$\{(a,b), (-b,-a), (a,\infty], [-\infty,-a) : 0<a<b\le \infty\}.$$
	Let $C_c^+(\overline{\R}_0)$ be the
	family of all non-negative continuous functions $g$ on $\overline{\R}_0$ with $g\equiv0$
	on $(-\delta,0)\cup(0,\delta)$ for some $\delta>0$.
	Denote by $\mathcal{M}(\overline{\R}_0) $ the space of all Radon
	measures on $\overline{\R}_0$
	endowed with the topology of vague convergence (denoted by $\overset{v}{\to}$),
	generated by the maps $\nu\to \int f d\nu$ for all $f\in C_c^+(\overline{\R}_0)$.
	For any $g\in\mathcal{B}_b^+(\overline{\R}_0)$, $\nu\in\mathcal{M}(\overline{\R}_0)$, we write $\nu(g):=\int_{\overline{\R}_0} g(x)\nu(\d x)$.
	A sequence of random elements $\nu_n$ in $\mathcal{M}(\overline{\R}_0)$ converges weakly to $\nu$,
	denoted as $\nu_n\overset{d}{\to}\nu, $
	if and only if for all $g\in C_c^+(\overline{\R}_0)$, $\nu_n(g)$ converges weakly to $\nu(g)$. 	Let $\mathcal{B}_1(\R)$ be the set of all the Borel functions $\varphi:\R\to[0,1]$
	with $\varphi\equiv1$ on $[-\delta, \delta]$  for some $\delta>0$.

	For any $x\neq0$ and   a measure $\nu$, we denote by $\nu/x$ the measure defined by
	$$(\nu/x)(g)=\int g(y/x)\nu(dy).$$

	It is well known that there exists
	a continuous  function $\tilde{L}: \R_+ \to \R_+$
	slowly varying at $\infty$ such that $\lim_{x\to\infty}\frac{\tilde{L}(x)}{L(x)}=1$,
	 $\tilde{L}(0+)\in (0, \infty)$
		 and $x^{-\alpha}\tilde{L}(x)$
	is strictly decreasing on $(0, \infty)$.
	In this paper, we always assume that
$L$ satisfies this  property.
Note that $\lim_{x\to 0} x^{-\alpha}{L}(x)=\infty$ and $\lim_{x\to \infty} x^{-\alpha}{L}(x)=0.$
Let $H(y):(0,\infty)\to(0,\infty)$
	be the inverse function of $x^{-\alpha}L(x)$.
Then
\begin{equation}\label{H-inverse}
H(y)^{-\alpha}L(H(y))=y, \quad y>0.
\end{equation}
It is well known (see, \cite[Theorem 1.5.12]{Bingham} for instance) that
		\begin{align}\label{formula-H}
		H(y)=y^{-1/\alpha}\bar{L}(y^{-1}),
	\end{align}
	with $\bar{L}$ being slowly varying at $\infty$.
	From now on, $\bar{L}$ always stands for the function above.
	
In \cite{YRR},
we studied the weak convergence of extremes of $\mathbb{X}$.
Let $h(t):=H(e^{-\lambda t})$, that is,
\begin{align}\label{def:h}
h(t)^{-\alpha}L(h(t))=e^{-\lambda t}.
\end{align}
Note that $h$ is strictly increasing.
Define
\begin{equation}\label{def-N}
	\mathcal{N}_t:=
	\mathbb{X}_t/h(t)=\sum_{v\in\mathcal{L}_t}\delta_{h(t)^{-1}\xi_t^v}.
\end{equation}	
In \cite[Theorem 1.1]{YRR},
we proved that $\mathcal{N}_t$
	converges weakly to a random measure $\mathcal{N}_\infty$. More precisely, for any $g\in C_c^+(\overline{\R}_0)$,
	\begin{align}\label{laplace-N}
		\lim_{t\to\infty}\E(e^{-\mathcal{N}_t(g)})&=E(e^{-\N_\infty(g)})=\E\left(\exp\left\{-C(e^{-g})W\right\}\right),
	\end{align}
where
\begin{align}\label{def:C}
C(\varphi):=\int_0^\infty e^{-\lambda r}\int_{\mathbb{R}_{0}}\E (1-\varphi(x)^{Z_{r}})v_\alpha(dx)\d r<\infty,\quad \forall \varphi\in\mathcal{B}_1(\R).
\end{align}
Moreover,
 $\mathcal{N}_\infty=\sum_{j}T_j\delta_{e_j},$
 where given $W$,
$\sum_{j}\delta_{e_j}$ is a Poisson random measure with intensity $\vartheta Wv_\alpha(dx)$,
$\{T_j,j\ge 1\}$
are  i.i.d. copies of a random variable $T$
with
\begin{align}\label{law-T}
P(T=k)
=\vartheta^{-1}\int_0^\infty  e^{-\lambda r}\P(Z_r=k)\d r,\quad k\ge1,
\end{align}
where
$v_\alpha(dx)$ is given by \eqref{def-v-alpha},
$\vartheta$ is defined in \eqref{d:lambda&vartheta},
and $\sum_{j}\delta_{e_j}$ and $\{T_j,j\ge 1\}$ are independent.
As a consequence, we proved in \cite[Corollary 1.2]{YRR}
that under $\P^*$, $\frac{R_t}{h(t)}
$ converges weakly. More precisely,
	\begin{align}\label{Rt/ht}
	\lim_{t\to\infty}\P^*\left(\frac{R_t}{h(t)}
	\le x\right)=\left\{\begin{array}{ll}
				 \E^*\left(e^{-\vartheta^* W x^{-\alpha}}\right), & x>0;\\
		0,&x\le 0.
	\end{array}\right.
	\end{align}

 Therefore, $\mathbb{X}_t$, normalized by $h(t)$, converges to a random measure. In particular,
 the largest position $R_t$ is of order $h(t)$ as $t\to\infty$.
 In this paper, we consider a function $\Lambda(t)$ which grows faster than
 the function
 $h(t)$ in the sense that $\lim\limits_{t \to \infty} \dfrac{\Lambda(t)}{h(t)}=\infty$, or slower than
 $h(t)$ in the sense that $\lim\limits_{t \to \infty} \dfrac{\Lambda(t)}{h(t)}=0$. When $\Lambda(t)$
 grows faster than $h(t)$, we
 find the  rate that $\P(R_t > \Lambda(t))$ converges to $0$  an $t\to\infty$, and describe the limit  of
 $\mathbb{X}_t/\Lambda(t)$,
 conditioned on $\{R_t > \Lambda(t)\}$,  at $t\to\infty$.  When $\Lambda(t)$ grows slower than
 $h(t)$, we
 find the rate that $\P(R_t \le \Lambda(t))$ converges to $0$ and describe the limit of  $\mathbb{X}_t/\Lambda(t)$  conditioned on $\{R_t\le \Lambda(t)\}$.
In this paper, we also study
the almost sure asymptotic behavior of
	$R_t$.
	 \bigskip
	
	\subsection{Main results}	\label{main results}
	
In this subsection, we state our main results. 	Let $\mathcal{H}(\R)$ denote the family of
uniformly continuous functions $\varphi:\R\to [0,1]$ with
$\varphi \equiv1$
in some neighborhood of $0$. Let $\mathcal{H}_0(\R)$ denote the family of all the  functions $\varphi\in\mathcal{H}(\R)$
with $\varphi \equiv 0$ on $(c, \infty)$  for some $c>0$.
	Note that if $g\in C_c^+(\bar{\R}_0)$ then $e^{-g}\in \mathcal{H}(\R)$.

\begin{theorem}\label{them:large-deviation}
If $\Lambda:[0,\infty)\to (0,\infty)$ satisfies
$\lim\limits_{t \to \infty} \dfrac{\Lambda(t)}{h(t)}=\infty$,
then for any $\varphi\in \mathcal{H}(\R)$,
\begin{align*}
	\lim\limits_{t\to \infty} e^{-\lambda t}\Lambda(t)^{\alpha}L(\Lambda(t))^{-1}\left( 1-\E\Big(\prod_{u\in\mathcal{L}_t}\varphi(\xi_t^u/\Lambda(t)\Big)\right)  =C(\varphi),
\end{align*}
where $C(\varphi)$ is defined in \eqref{def:C}.
In particular,
\begin{align}
	\lim\limits_{t\to \infty} e^{-\lambda t}\Lambda(t)^{\alpha}L(\Lambda(t))^{-1} \P(R_t > \Lambda(t))=\vartheta^*,
\end{align}
where
$\vartheta^*$ is defined in \eqref{d:varthetastar}.
\end{theorem}

In \cite{Shiozawa2}, Shiozawa
studied the upper deviation of  the maximal displacement of a branching symmetric stable process with spatially inhomogeneous branching structure, and proved some weak convergence results.
	
\begin{corollary}\label{cor:large-deviation}
	If $\Lambda:[0,\infty)\to (0,\infty)$ satisfies
	$\lim\limits_{t \to \infty} \dfrac{\Lambda(t)}{h(t)}=\infty$,
then conditioned on $\{R_t>\Lambda(t)\}$,
\begin{itemize}
	\item [(1)]$R_t/\Lambda(t)$ converges weakly to a random variable $R^*$ with density
	$\alpha x^{-1-\alpha}{\bf 1}_{(1,\infty)}(x).$
	\item[(2)] $\mathbb{X}_t/\Lambda(t)$ converges  weakly to $T\delta_{R^*}$,
where the law of $T$ is given in \eqref{law-T},
 and $T$  and $R^*$ are  independent.
\end{itemize}
\end{corollary}

Recall that  $f(s)=\sum_{k}p_k s^k$
and $q=\P(\mathcal{S}^c)$. Note that $f'(q)\in[0,1)$. Put
\begin{align}\label{def:rho}
	 \rho:=\beta(1-f'(q)).
\end{align}
For any $\theta\ge0$, define
\begin{equation}\label{def-phi}
\phi(\theta)=\E(e^{-\theta W}).
\end{equation}

\begin{theorem}\label{them:lower deviation}
	Let  $\Lambda:[0,\infty)\to (0,\infty)$ be  a non-decreasing function. Assume that
	$$\sum_{n=1}^\infty n \Lambda(n)^{-\alpha}L(\Lambda(n))<\infty \mbox{ and }\lim_{t\to\infty}\frac{\Lambda(t)}{h(t)} =0,
	$$	
and if $p_k=0$ for all $k\ge3$, we further assume that $\Lambda(t)>e^{\gamma t}$
with some $\gamma>0$ for $t$ sufficiently large.
	Then for any $\varphi\in\mathcal{H}_0(\R)$,
	$$\lim_{t\to\infty}e^{\rho (t-r(t))}\E^*\Big(\prod_{u\in\mathcal{L}_t}\varphi(\xi_t^u/\Lambda(t))\Big) = \frac{1}{1-q} A\left[\phi(C(\varphi))\right] ,$$
	where $r(t)$ is defined by $h(r(t))=\Lambda(t)$, $C(\varphi)$ is defined in \eqref{def:C}, and $A(s)$ is defined in \eqref{limit-F} below.
	In particular,
	$$\lim_{t\to\infty}e^{\rho (t-r(t))}
 \P^*\Big(R_t\le \Lambda(t)\Big)
 = \frac{1}{1-q} A\left[\phi(\vartheta^*)\right] .$$
\end{theorem}

\begin{corollary}\label{cor:lower deviation}
	Let  $\Lambda:[0,\infty)\to (0,\infty)$ be  a non-decreasing function. If
	$$\sum_{n=1}^\infty n \Lambda(n)^{-\alpha}L(\Lambda(n))<\infty \mbox{ and }\lim_{t\to\infty}\frac{\Lambda(t)}{h(t)}=0,$$  then under $\P^*$, conditioned on $\{R_t\le \Lambda(t)\}$, $\mathbb{X}_t/\Lambda(t)$ converges weakly to some random measure
	 $\Xi=\sum_{k=1}^\mathcal{K} \bar{\mathcal{N}}^{(k)}_\infty,$ where
	\begin{itemize}
		 \item[(i)] $\{\bar{\mathcal{N}}^{(k)}_\infty, k\ge 1\}$ are i.i.d. with the same law as
    $P(\mathcal{N}_{\infty}\in\cdot|\mathcal{N}_\infty(\R)\neq 0, \mathcal{N}_{\infty}((1,\infty))=0);$
		 \item[(ii)] $\mathcal{K}$ is a positive integer valued random variable with generating function $$E\left(s^\mathcal{K}\right)=\frac{A( (\phi(\vartheta^*)-q)s+q)}{A(\phi(\vartheta^*))};$$
      \item[(iii)] $\{\bar{\mathcal{N}}^{(k)}_\infty, k\ge 1\}$  and $\mathcal{K}$ are independent.
\end{itemize}
\end{corollary}

\begin{remark}
(1) It is interesting that, by Theorem \ref{them:large-deviation} and Corollary \ref{cor:large-deviation}, conditioned on $\{R_t>\Lambda(t)\}$, the limits of $R_t/\Lambda(t)$ and $\mathbb{X}_t/\Lambda(t)$ do not depend on the function $\Lambda(t)$: for any $\Lambda(t)$ satisfying $\lim\limits_{t \to \infty} \dfrac{\Lambda(t)}{h(t)}=\infty$, the limits are the same. The limit of the point process $\mathbb{X}_t/\Lambda(t)$, conditioned on $\{R_t>\Lambda(t)\}$, is a point measure supported on one point $R^*$.

(2) By Theorem \ref{them:lower deviation} and and Corollary \ref{cor:lower deviation}, the limit of $\mathbb{X}_t/\Lambda(t)$, conditioned on the event $\{R_t\le \Lambda(t)\}$, does not depend on $\Lambda(t)$: for any $\Lambda(t)$ satisfying $\lim\limits_{t \to \infty} \dfrac{\Lambda(t)}{h(t)}=0$ and $\sum_{n=1}^\infty n \Lambda(n)^{-\alpha}L(\Lambda(n))<\infty$, the limit $\Xi$ does not depends on $\Lambda$.
Comparing $\Xi$  with $\mathcal{N}_{\infty}$ (the limit of  $\mathbb{X}_t/h(t)$),
we see that $\Xi$ is a random
sum of independent copies of $\bar{\mathcal{N}}^{(k)}_\infty$, with common law equal to that of
$\mathcal{N}_{\infty}$ condition on $\{\mathcal{N}_{\infty}(\R)\neq 0, \mathcal{N}_{\infty}((1,\infty))=0\}$.

(3)
In the special case that   $L\equiv 1$, we have $h(t)=e^{\frac{\lambda}{\alpha} t}$.
Consider $\Lambda(t)=e^{c\frac{\lambda}{\alpha} t}$ for some constant  $c>0$. If $c>1$,
by Theorem \ref{them:large-deviation},
	$$\lim\limits_{t\to \infty} e^{(c-1)\lambda t} \P(R_t > \Lambda(t))=\vartheta^*.$$
	If $0<c<1$, then $r(t)=ct$.
By Theorem \ref{them:lower deviation},
	$$\lim_{t\to\infty}e^{ ((1-c)\rho t}\P^*(R_t\le \Lambda(t)) =A\left(\phi\left(\vartheta^* \right)  \right).$$
\end{remark}

We now state some almost sure convergence results of $R_t$.

	\begin{theorem}\label{main-theorem-inf} It holds that
		\begin{align*}
			\liminf\limits_{t \to \infty}\dfrac{R_t }{H(e^{-\lambda t}\log t)} = \left(\vartheta^* W\right)^{\frac{1}{\alpha}}, \qquad \P^* \textup{-a.s.}
		\end{align*}
	\end{theorem}

	\begin{theorem}\label{them:upper-limit}
	Suppose that  $G:[0,\infty)\to (0,\infty)$ is a non-decreasing function
satisfying
$\lim_{t\to\infty}\frac{G(t)}{h(t)}=\infty$.
	\begin{itemize}
		\item[(1)] If $\sum_{n}e^{\lambda n}G(n)^{-\alpha}L(G(n))<\infty$, then
		\begin{align*}
			\limsup\limits_{t \to \infty}\dfrac{ R_t}{G(t)}= 0,  \qquad \P^* \textup{-a.s.}
		\end{align*}
		\item[(2)] If  $\sum_{n}e^{\lambda n}G(n)^{-\alpha}L(G(n))=\infty$, then
		\begin{align*}
			\limsup\limits_{t \to \infty}\dfrac{ R_t}{G(t)}= \infty,  \qquad \P^* \textup{-a.s.}
		\end{align*}	
	\end{itemize}
\end{theorem}

This implies that, for a non-decreasing function $G$ satisfying
$\lim_{t\to\infty}\frac{G(t)}{h(t)}=\infty$, either $\limsup\limits_{t \to \infty}\dfrac{ R_t}{G(t)}= 0$ $\P^* \textup{-a.s.}$, or $\limsup\limits_{t \to \infty}\dfrac{ R_t}{G(t)}= \infty$ $\P^* \textup{-a.s.}$.
Similar result holds for a  subordinator. If $\xi$ is a subordinator with infinite mean, then for any $G:[0,\infty)\to (0,\infty)$ being an increasing function such that $\frac{G(t)}{t}$ increases, then $\limsup_{t\to\infty}\frac{\xi_t}{G(t)}= 0 \mbox{ or }\infty$ almost surely, see \cite[Theorem 13]{Bertoin}.
	
\bigskip

Now we give some intuitive idea for
one of the main techniques of this paper.
	Let
	    $Y_j:=\xi_{j}-\xi_{j-1}, j\ge1$.
	Then $\{Y_j\}$ are
	i.i.d. It is easy to see from \eqref{24}  that
	$$\mP(\xi_n>x)\sim n
	\mP(\xi_1>x)\sim\mP(\max_{j\le n}Y_j>x),\quad x\to\infty.
	$$
Thus the maximum $\max_{1\le j\le n}Y_j$ plays
a dominating role in the asymptotic behavior of $\xi_n$.

For $t\ge 0$, $u\in \mathbb{T}$,  we set $X_{u, t}:=\xi^u_{\sigma_u\wedge t}-\xi^u_{b_u\wedge t}$.  Then  we have
$$\mathbb{X}_t=\sum_{u\in\mathcal{L}_t} \delta_{\xi_t^u}=\sum_{u\in\mathcal{L}_t}\delta_{\sum_{v\in I_u}X_{v,t}}.$$
We will see in
Lemma \ref{lem:one-big'}
that the asymptotic behavior of $\mathbb{X}_t$ is governed  by $\mathbb{Y}_t$:
\begin{align*}
\mathbb{Y}_t :=\sum_{u\in\mathcal{L}_t}\sum_{v\in I_u}\delta_{X_{v,t}}.
\end{align*}
Thus, to prove Theorems \ref{them:large-deviation},    \ref{main-theorem-inf} and \ref{them:upper-limit},  we first establish the corresponding results for $\mathbb{Y}_t$. This technique  has been employed in \cite{YRR} for branching L\'evy processes and
in \cite{BHR, BHR2,Durrett83}
for branching random walks with heavy tails.  However this technique (Lemma \ref{lem:one-big'})) does not work for the proof of lower deviation result in Theorem \ref{them:lower deviation}.  We establish
lower large deviation results of  $\mathbb{X}_t$ and $\mathbb{Y}_t$ separately,
and it turns out that the results are identical,
which is somewhat surprising.

For  branching random walks,  several authors have studied
the convergence of  the extremes under an exponential
moment assumption on the  displacements of the offspring from the parent,
see A\"{i}d\'{e}kon \cite{Aldekon},
Hu and Shi \cite{HS},
and Madaule \cite{Madaule}. Recently, many researchers
studied related topics for
branching random walks with heavy-tailed displacements. Assume that the displacements of the offspring from the parents are i.i.d. with
$$
P(X>x)\sim ae^{-L(x)x^{r}},
\quad x\to\infty,
$$
where $a>0$, $L$ is slowly varying at $\infty$ and $r\in[0,1)$.
When $r\in(0,1)$, the maximum $M_n$ grows polynomially. For example,
if $L$ is a constant,
then $M_n/n^{1/r}$ converges to a positive constant almost surely.
See, \cite{DGT23,DG,DGT20,Gantert} for
more related results.
When $r=0$,  $\log t/L(t)\to0$ (or $L(t)=o(\log t)$)
as $t\to\infty$,
 the extremes have been investigated in \cite{BDGP}.
When $r=0$ and $L(x)=\alpha\log x-\log \tilde{L}(x)$ where $\tilde{L}$ is slowly varying at $\infty$,
Durrett \cite{Durrett83}
proved that $a_n^{-1}M_n$  converges weakly,  where $a_n=m^{n/\alpha}L_0(m^n)$ and $L_0$ is slowly varying at $\infty$. Recently,  the extremal processes of the branching random walks with regularly varying steps were studied by Bhattacharya et al. \cite{BHR,BHR2}.
It was proved in \cite{BHR,BHR2}
that the point random measures $\sum_{|v|=n}\delta_{a_n^{-1}S_v}$, where $S_v$ is the position of $v$, converges weakly to a Cox cluster process, which are quite different from the case with exponential moments.
Recently, Bhattacharya \cite{B} studied the large deviations of extremes in branching  random walk with regularly varying  displacements, corresponding to our results for $\Lambda(t)$ growing fast than $h(t)$.

\section{Upper deviation of  $\mathbb{X}_t$ and $\mathbb{Y}_t$ }

	It is well known (see \cite[Theorem 1.5.6]{Bingham} for instance) that, for any $\epsilon>0$, there
	exists $a_\epsilon>0$ such that for any
	$x, y>a_{\epsilon}$,
	\begin{align}\label{22}
		\frac{L(y)}{L(x)}\le 2\max\{(y/x)^{\epsilon}, (y/x)^{-\epsilon}\}, \qquad \frac{\bar{L}(y)}{\bar{L}(x)}\le 2\max\{(y/x)^{\epsilon} , (y/x)^{-\epsilon}\}.
	\end{align}

	Let  $C_{b}^0(\R)$ be the space of all  bounded continuous functions vanishing in a neighborhood of $0$.  Recall the definition of $v_\alpha$ in \eqref{def-v-alpha}. In the following lemma, we present a generalization of \cite[Lemma 2.1]{YRR}. Since the proof follows a similar line of reasoning, we omit it.
	
	\begin{lemma}\label{lem:vague}
If  {\bf Assumption \ref{assum2} }holds,
then  for any $g\in C_b^0(\R)$ and $s>0$,
		\begin{equation}\label{vague-conv1}
			\lim_{x\to\infty} x^{\alpha}L(x)^{-1}\mE\left(g\left(\frac{\xi_s}{x}\right)\right)=s\int_{\R_0} g(x)v_\alpha(\d x).
		\end{equation}
	\end{lemma}

	In  \cite[Lemma 2.2]{YRR}, we have proved that, under  {\bf Assumption 2}, there exist $c_0>0$ and $x_0>0$ such that for any $s>0$ and $x>x_0$,
	\begin{align}\label{tail-prob}
			\mP(|\xi_s|>x)\le c_0 s x^{-\alpha}L(x).
	\end{align}
Using \eqref{22} with $\epsilon=1$, we see that,
for any $c\in(0,1)$ and $y>0$
sufficiently small so that
$cH(y)>a_1+x_0$, it holds that
\begin{align}\label{tail-prob2}
	\mP(|\xi_s|>cH(y))\le
	&c_0 s c^{-\alpha}H(y)^{-\alpha}L(cH(y))\nonumber\\
\le& 2c_0 s c^{-\alpha-1}H(y)^{-\alpha}L(H(y))=2c_0  c^{-\alpha-1}sy,
\end{align}
where in the last equality we used \eqref{H-inverse}.
It was shown in \cite{Eric} that, for any $t>0$,
\begin{align}\label{sup-xi}
\lim_{x\to\infty}x^{\alpha}L(x)^{-1}\mP\left(\sup_{0\le s\le t}\xi_s>x\right)=\lim_{x\to\infty}x^{\alpha}L(x)^{-1} \mP(\xi_t>x)=\frac{q_1}{\alpha}
\end{align}
and
\begin{align}\label{inf-xi}
	\lim_{x\to\infty}x^{\alpha}L(x)^{-1}\mP\left(\inf_{0\le s\le t}\xi_s<-x\right)=\lim_{x\to\infty}x^{\alpha}L(x)^{-1} \mP(\xi_t<-x)=\frac{q_2}{\alpha}.
\end{align}

 We now recall a special many-to-one formula. For more general many-to-one formulas, see \cite[Theorem 8.5]{Hardy-Harris}.
 For any $u\in\mathbb{T}$,  let $n^u$ be the number  of particles in $I_u\setminus \{o\}$.
	
	\begin{lemma}[Many-to-one formula]\label{many-to-one}
		Let $\{n_t\}$ be a Poisson process with parameter $\beta$ on some probability space $(\Omega, {\cal G}, P)$.
		Then for any $g\in\mathcal{B}_b^+(\R)$,
		$$\E\left(\sum_{v\in\mathcal{L}_t} g(n^v)\right)
		=e^{\lambda t}
		E\left(g(n_t)\right).
		$$
	\end{lemma}

	\subsection{A key lemma on ``one big jump"}

	In the remainder of the paper, we use $g_1(t)\overset{t}{\le}g_2(t) $ to denote that $g_1(t)\le g_2(t)$  for sufficiently large $t$. For any  nonnegative  function $g$ and measure $\nu=\sum_{k=1}^n \delta_{x_k}$, define
	$$
	\mathcal{I}(g,\nu):=\prod_{k=1}^n g(x_k).
	$$
	Here we use the convention that $\prod_{k=1}^0 g(x_k)=1$.
	It is clear that $\cI(g,\nu)=e^{\nu(\log g)}.$

	Recall that, for $t\ge 0$ and $u\in \mathbb{T}$,
	$$X_{u, t}=\xi^u_{\sigma_u\wedge t}-\xi^u_{b_u\wedge t}.$$
		Let $\mathcal{D}_t:=\{u\in\mathbb{T}:b_u\le t, Z_t^u>0\}$, where $Z_t^u$ is the number of  offspring of particle $u$ alive at time $t$.
Define
\begin{align} \label{M_t}
	M_t :=\max_{v\in\mathcal{L}_t}\max_{u\in I_v}X_{u,t}=\max\limits_{u\in\mathcal{D}_t}X_{u,t}
\end{align}
Recall that
\begin{align}
	\mathbb{Y}_t=\sum_{v\in \mathcal{L}_t}\sum_{u\in I_v} \delta_{X_{u,t}}=\sum_{u\in\mathcal{D}_t}Z^u_t \delta_{X_{u,t}}.
\end{align}
Let $0<s<t$.
The particles in $\mathcal{D}_t$ can be divided into two groups: those born before time $t-s$ and those born after $t-s$.
We define
\begin{align}\label{def:M1}
	M_{s,t}:=\max\limits_{u \in \mathbb{T}:b_u\leq t-s}|X_{u,t}|
\end{align}	
and
\quad
\begin{align}\label{def:Y-st}
	\mathbb{Y}_{s,t}:=\sum_{u\in\mathcal{D}_t:t-s< b_u\le t}Z^u_t \delta_{X_{u,t}}=\sum_{u\in\mathbb{T}:t-s< b_u\le t}Z^u_t \delta_{X_{u,t}}.
\end{align}
Using the tree structure, we can categorize all particles born after $t-s$ according to the branches formed by particles that were alive at time
$t-s$. More precisely,
\begin{align} \label{314}
	\mathbb{Y}_{s,t}=\sum_{v\in\mathcal{L}_{t-s}}\sum_{u:v<u, b_u\le t}Z^u_t \delta_{X_{u,t}}=:\sum_{v\in\mathcal{L}_{t-s}}\mathbb{Y}^v_{s,t}.
\end{align}	
By the branching property and the Markov property, conditioned on $\mathcal{F}_{t-s}$, $\{\mathbb{Y}^v_{s,t}, v\in\mathcal{L}_{t-s}\}$ are i.i.d. with a common law
equal to that of
\begin{align}\label{def:Y-stprime}
	\mathbb{Y}'_s:=\sum_{u\in \mathcal{D}'_s}Z^u_s \delta_{X_{u,s}},
\end{align}	
where $\mathcal{D}'_s:=\mathcal{D}_s \setminus \{o\}.$
We will also use the following notation later:
\begin{align} \label{M_t'}
	M_t' :=\max\limits_{u\in\mathcal{D}'_t}X_{u,t}.
\end{align}

\begin{lemma}\label{lem:m1}
If $a(t)$ is  a positive function with  $\lim\limits_{t \to \infty} a(t)=\infty$, then
	\begin{align}
		\lim\limits_{s \to \infty} \limsup\limits_{t \to \infty} e^{-\lambda t}a(t)^{\alpha}L(a(t))^{-1} \mathbb{P}(M_{s,t}>a(t) )=0.
	\end{align}
\end{lemma}	
\noindent {\bf Proof:}
When $a(t)=h(t)$,
it has been proven in \cite[(2.14)]{YRR}) that
\begin{align}\label{2.4.1}
	\mathbb{P}(M_{s,t}>a(t) )\overset{t}{\le }c_0e^{\lambda t}a(t)^{-\alpha}L(a(t))\Big(\lambda^{-1}e^{-\lambda s}+\frac{e^{-\beta s}-e^{-\lambda s}}{\lambda-\beta}\Big).
\end{align}
For the general case, \eqref{2.4.1} also holds. The proof follows almost the same argument, so we omit the details. Then the desired result follows immediately.

\hfill$\Box$

	The following key lemma says that $\mathbb{X}_t$ and $\mathbb{Y}_t$ have similar  asymptotic behaviors.
In the proof of Lemma \ref{lem:one-big'}, we show that with high probability, for all  $v\in\mathcal{L}_t$, there exists at most one particle  $u\in I_v$ that experiences a ``large jump".

\begin{lemma}\label{lem:one-big'}
If $a(t)$ is a positive function such that $a(t)\overset{t}{>}e^{\epsilon t  }$ for  some $\epsilon\in (0, 1)$, then
for any $\varphi\in \mathcal{H}(\R)$,
	\begin{align}\label{901}
		\lim_{t\to\infty}e^{-\lambda t} a(t)^{\alpha}L(a(t))^{-1}\E\left( \left|\mathcal{I}(\varphi,\mathbb{X}_t/a(t))-\mathcal{I}( \varphi,\mathbb{Y}_t/a(t))\right|\right) =0.
	\end{align}
	
\end{lemma}
	\noindent\begin{proof}
		We divide the proof into three steps.
		
		\textbf{Step 1}
		For any $t>0$ and $\theta\in(0,1]$,
				define
	 $$
	 A_t(\theta):=
	 \displaystyle\bigcap\limits_{v \in \mathcal{L}_t}\textstyle\left\{\sum\limits_{u \in I_v}\mathbf{1}_{\{|X_{u,t}| >\theta a(t) /\log a(t)\}} \leq 1\right\}.$$
	     We claim that  for any $p\in(0,2)$,
		\begin{align} \label{31-1}
			\lim_{t\to\infty}e^{-\lambda t} a(t)^{p\alpha} \mathbb{P}(A_t(\theta)^c) = 0.
		\end{align}
		\noindent Note that
		\begin{align} \label{32}
			\P(A_t(\theta)^c|\mathcal{F}_t^{\mathbb{T}})  \leq \sum_{v \in \mathcal{L}_t}
			\P\left(\left\{\sum\limits_{u \in I_v}\mathbf{1}_{\{|X_{u,t}| > \theta a(t)  /\log a(t)\}} \geq 2\right\}|\mathcal{F}_t^{\mathbb{T}}\right).
		\end{align}
By \eqref{tail-prob} we have that
		\begin{align} \label{33}
			&\P\left(|X_{u,t}| > \theta a(t) /\log a(t) | \mathcal{F}_t^{\mathbb{T}}\right) =\mP(|\xi_s|>\theta a(t) /\log a(t))|_{s=\tau_{u,t}}\nonumber\\
			 \overset{t}{\le } &\left(c_0  \theta^{-\alpha}t\cdot a(t)^{-\alpha} [ \log a(t)]^{\alpha} L(\theta a(t)/\log a(t) )\right)\wedge 1=:p_t.
		\end{align}
	Recall that,  for any $v\in\mathbb{T}$,   $n^v$ is the number  of particles in $I_v\setminus \{o\}$. Thus $|I_v|=n^v+1$.
Since,  conditioned on $\FF_t^{\T}$ , $\{X_{u,t},u\in I_v\}$ are independent, we have
		\begin{align*}
			&\mathbb{P}\left(\sum_{u \in I_v} \mathbf{1}_{\{|X_{u,t}| > \theta a(t) /\log a(t)\}}\geq 2|\FF_t^{\T}\right)  \leq \sum_{m=2}^{n^v+1} \binom{n^v+1}{m} p_t^m \\
			\leq &p_t^2 \sum_{m=0}^{n^v-1} n^v (n^v+1) \binom{n^v-1}{m} p_t^{m} \\
			=& p_t^2 n^v (n^v+1)(1+p_t)^{n^v-1}.
		\end{align*}
		\noindent Thus by \eqref{32} and the many-to-one formula, we have
		\begin{align} \label{34}
			\P(A_t(\theta)^c) &
			\leq  e^{\lambda t} p_t^2 E(n_t(n_t+1)(1+p_t)^{n_t-1})  =e^{\lambda t}  p_t^2(2+(1+p_t)\beta t)\beta te^{\beta t p_t}
\nonumber\\ &\sim (\beta t)^2e^{\lambda t}  p_t^2, \quad t\to\infty.
		\end{align}
	Here we used the fact $tp_t\to0$ as $t\to\infty$.
	Now  \eqref{31-1} follows from \eqref{33} and \eqref{34} immediately.
		
		\textbf{Step 2}
		In the remainder of this proof, we fix a constant $c>2\alpha+e^2\beta/\epsilon$.
		Define
		$B_t :=\displaystyle\bigcap\limits_{v \in \mathcal{L}_t} \{n^v \leq  c\log a(t)\}$. Using the many-to-one formula, we get
		\begin{align}\label{3.1}
			\P(B_t^c) & =  \P\left(\bigcup\limits_{v \in \LL_t} \{n^v > c \log a(t)\}\right) \leq \E\left(\sum_{v \in \LL_t} \mathbf{1}_{\{n^v >c  \log a(t)\}}\right) \notag\\
			&= e^{\lambda t}P(n_t> c\log a(t))\le e^{\lambda t}\inf\limits_{r>0} e^{-r c\log a(t)} E(e^{rn_t}) \notag\\
			& =  e^{\lambda t} e^{-(\text{log} \log a(t)+\log c-\text{log}\beta t-1) c\log a(t)-\beta t}.
		\end{align}
		Since $a(t)\overset{t}{>} e^{\epsilon t}$, we have
		 $$\log \log a(t)+\log c\overset{t}{>}  \log( c\epsilon t)\ge  \log(\beta t)+2. $$
		 Thus by \eqref{3.1} , we have
		$$\P(B_t^c) \overset{t}{<}  e^{\lambda t} e^{-c\log a(t)}e^{-\beta t}=e^{\lambda t}a(t)^{-c} e^{-\beta t}.$$
		Since $c>2\alpha$, we get
		\begin{align}\label{Bto0}
			\lim_{t\to\infty}e^{-\lambda t}a(t)^{2\alpha}\P(B_t^c)=0.
		\end{align}

		\textbf{Step 3}
		Let $v^{\prime} \in I_v$ be such that $|X_{v^{\prime},t}| = \max\limits_{u \in I_v}\{|X_{u,t}|\}$. We note that, on the event $A_t(\theta )$, $|X_{u,t}| \leq  \theta a(t) /\log a(t) $ for any $u \in I_v \setminus \{v^{\prime}\}$.
		Since $\xi_t^v = \sum_{u \in I_v} X_{u,t}$, on the event
		 $A_t (\theta)\cap B_t$,
		we have that
		\begin{align} \label{35}
			\left|\xi_t^v-X_{v^{\prime},t}\right|=\left|\sum_{u \in I_v \setminus \{v^{\prime}\}}X_{u,t}\right| \leq
			\frac{n^v \theta a(t)}{ \log a(t)}
			\leq c\theta  a(t).
		\end{align}

	For  $\varphi \in \mathcal{H}(\R)$, we have
		$\varphi \equiv1$ on $[-\delta, \delta]$ for some $\delta> 0$.
	Since $\varphi$ is uniformly continuous,
	for any $\gamma>0$,
	there exists $\eta>0$ such that
	$|\varphi(x_1)-\varphi(x_2)|\le \gamma$
	whenever $|x_1-x_2|<\eta$. We now fix an arbitrary $\gamma>0$ and
		the corresponding $\eta$.
	
		Recall $c>2\alpha+e^2\beta/\epsilon$. Choose
	$\theta$  small enough so that  $c\theta<\eta\wedge(\delta/2)$. We assume that $t$ is sufficiently large so that $\log a(t)>2\theta /\delta$.
	We note that, on the event $A_t(\theta)$, $|X_{u,t}|\le \theta a(t)/\log a(t)\le a(t)\delta/2$ for any
	$u\in I_v\setminus\{v'\}$,
	and thus  $\varphi(X_{u,t}/a(t))=1$. It follows that on the event
	$A_t (\theta)\cap B_t$,
	\begin{align}\label{903}
		\left|\prod_{v\in\mathcal{L}_t}\varphi(\xi_t^v/a(t))- \prod_{v\in\mathcal{L}_t}
		\prod_{u\in I_v}
		\varphi(X_{u,t}/a(t))\right|&=\left|\prod_{v\in\mathcal{L}_t}\varphi(\xi_t^v/a(t))- \prod_{v\in\mathcal{L}_t}\varphi(X_{v',t}/a(t))\right|\nonumber\\
		&\le \sum_{v\in\mathcal{L}_t} |\varphi(\xi_t^v/a(t))-\varphi(X_{v',t}/a(t))|.
	\end{align}
By \eqref{35}, on $A_t(\theta)\cap B_t$,  we have
	\begin{align} \label{35'}
	\left|\xi_t^v-X_{v^{\prime},t}\right|/a(t)
	\leq c\theta\le  \eta\wedge(\delta/2).
\end{align}
Thus if $|X_{v',t}|\le \delta a(t)/2$, then  $|\xi_t^v|/a(t)< \delta$, which implies that $\varphi(\xi^v_t/a(t))-\varphi(X_{v',t}/a(t))=0$.
Hence  by \eqref{903},
on $A_t(\theta)\cap B_t$,
	\begin{align}
		&\left|\mathcal{I}(\varphi,\mathbb{X}_t/a(t))-\mathcal{I}(\varphi,\mathbb{Y}_t/a(t))\right| \le
		\gamma
		\sum_{v\in\mathcal{L}_t} {\bf 1}_{\{|X_{v',t}|>a(t)\delta/2\}}\\
		\le&\gamma
		\sum_{v\in\mathcal{L}_t}\sum_{u\in I_v}{\bf 1}_{\{|X_{u,t}|>a(t)\delta/2\}}
		=\gamma\sum_{u\in\mathcal{D}_t}Z^u_tg_t(X_{u,t}),
	\end{align}	
where $g_t(y)={\bf 1}_{|y|>a(t)\delta/2}$.
Therefore, for any fixed
$s\in (0, t)$, on $A_t(\theta)\cap B_t\cap \{M_{s, t}\le a(t)\delta/2\}$, we have
$$
\left|\mathcal{I}(\varphi,\mathbb{X}_t/a(t))-\mathcal{I}(\varphi,\mathbb{Y}_t/a(t))\right|
\le \gamma\sum_{u:t-s<b_u\le t}Z^u_tg_t(X_{u,t})=\gamma\mathbb{Y}_{s,t}(g_t).$$
	Since $\varphi$ takes values in $[0,1]$,  we have
	\begin{align}\label{2.5.1}
		&\E\left|\mathcal{I}(\varphi,\mathbb{X}_t/a(t))-\mathcal{I}(\varphi,\mathbb{Y}_t/a(t))\right| \nonumber\\
		\le &\P(A_t(\theta)^c)+\P(B_t^c)+\P(M_{s,t}>a(t)\delta/2)
		+\gamma
		\E(\mathbb{Y}_{s,t}(g_t)).
	\end{align}
	By \eqref{314} and \eqref{tail-prob}, we have that
	\begin{align}
		\E(\mathbb{Y}_{s,t}(g_t))=&\E(Z_{t-s})\E(\mathbb{Y}'_s(g_t))=e^{\lambda(t-s)}\E\left(\sum_{u\in\mathcal{D}_s'}Z^u_s\P(|X_{u,s}|>a(t)\delta/2|\mathcal{F}_s^{\mathbb{T}})\right)\\
		\overset{t}{\le}&c_0(\delta/2)^{-\alpha}a(t)^{-\alpha}L(\delta a(t)/2)e^{\lambda(t-s)}\E\left(\sum_{u\in\mathcal{D}_s'}Z^u_s\tau_{u,s}\right)\\
		\le &c_0(\delta/2)^{-\alpha}a(t)^{-\alpha}L(\delta a(t)/2)e^{\lambda t}s,
	\end{align}
where in the last inequality we used the following inequality:
	$$\E\left(\sum_{u\in\mathcal{D}_s'}Z^u_s\tau_{u,s}\right)\le \E\left(\sum_{u:b_u\le s}Z^u_s\tau_{u,s}\right)\le se^{\lambda s},$$
which follows from the display below (2.19) in the proof of  \cite[Proposition 2.1]{YRR}.
	Thus combining \eqref{2.5.1}, \eqref{31-1} and \eqref{Bto0}, we have
	\begin{align}
		&\limsup_{t\to\infty} e^{-\lambda t} a(t)^{\alpha}L(a(t))^{-1}\E\left|\mathcal{I}(\varphi,\mathbb{X}_t/a(t))-\mathcal{I}(\varphi,\mathbb{Y}_t/a(t))\right| \\
		\le& \limsup\limits_{t \to \infty} e^{-\lambda t}a(t)^{\alpha}L(a(t))^{-1} \mathbb{P}(M_{s,t}>a(t) \delta/2)+c_0(\delta/2)^{-\alpha}s\gamma.
	\end{align}
Letting $\gamma\to0$ first, and then letting $s\to\infty$ and applying Lemma \ref{lem:m1}, we arrive at the desired result.

	\hfill$\Box$
	\end{proof}

	\subsection{Proof of Theorem \ref{them:large-deviation}}\label{prop 3.1}
	
	We emphasize here that the definition of $M_{s,t}$ in \cite[Section 2.3]{YRR} coincides with $\mathbb{Y}'_s/h(t)$.

	\begin{lemma}\label{lem:Ms}
	For any $s>0$ and $\varphi\in\mathcal{H}(\R)$,
	\begin{align}\label{251}
	\lim_{x\to\infty}x^\alpha L(x)^{-1}\left[ 1-\E
	\left(\mathcal{I}( \varphi,\mathbb{Y}'_s/x)\right)\right] =\int_0^s \left(e^{\lambda (s-r)}-e^{-\beta (s-r)}\right)\int_{\R}\E \left(1-\varphi(y)^{Z_{r}}\right)v_\alpha(dy)\,\d r.
	\end{align}
	\end{lemma}
\noindent{\bf Proof:} It has been proven in \cite[(2.17) and the first display on page 636]{YRR} that for $\varphi(x)=e^{-g(x)}$ with $g\in C_0^+(\bar{\R}_0)$,
\begin{align*}
	&\lim_{t\to\infty}h(t)^{\alpha}L(h(t))^{-1}\left[ 1-\E
	\left(\mathcal{I}( \varphi,\mathbb{Y}'_s/h(t))\right)\right]\\
	 &=\int_0^s \left(e^{\lambda (s-r)}-e^{-\beta (s-r)}\right)\int_{\R}\E \left(1-\varphi(y)^{Z_{r}}\right)v_\alpha(dy)\,
	  \d r.
	\end{align*}
By examining the proof of the above limit, we observe that it holds for any $\varphi\in\mathcal{H}(\R)$.
The desired result now follows.

\hfill$\Box$

	Suppose $\varphi \in\mathcal{B}_1(\R)$. If $\{\varphi_n, n\ge 1\}\subset \mathcal{B}_1(\R)$
are such that 	$\varphi$ and $\varphi_n$ are identically 1 on $[-\delta, \delta]$ for some $\delta>0$, and that $\varphi_n\to \varphi$ almost everywhere,
	then by the dominated convergence theorem,
\begin{align}
	C(\varphi_n)\to C(\varphi).
\end{align}
Hence $C(\varphi)$ is continuous in $\varphi$.
Recall $M_t$ is defined in \eqref{M_t}.

\begin{proposition}\label{prop:large-M}
If $\Lambda(t)$ is a positive function with  $\lim\limits_{t \to \infty} \dfrac{\Lambda(t)}{h(t)}=\infty$, then
for any $\varphi\in \mathcal{H}(\R)$,
	\begin{align}\label{lim-Yt}
		\lim\limits_{t\to \infty}	e^{-\lambda t}\Lambda(t)^{\alpha}L(\Lambda(t))^{-1}\left(1-\E\left(\mathcal{I}(\varphi, \mathbb{Y}_t/\Lambda(t))\right) \right)=C(\varphi),
	\end{align}
	where $C(\varphi)$ is defined in \eqref{def:C}.
	Furthermore,
	\begin{align}\label{lim-Mt}
	\lim\limits_{t\to \infty}	e^{-\lambda t}\Lambda(t)^{\alpha}L(\Lambda(t))^{-1}\P(M_t>\Lambda(t))
	=\vartheta^*,
	\end{align}
	where $\vartheta^*$ is defined in \eqref{d:varthetastar}.
\end{proposition}	
	
\noindent	{\bf Proof :}
For  $\varphi \in \mathcal{H}(\R)$,
we have $\varphi\equiv1$ on $[-\delta, \delta]$ for  some $\delta> 0$.
It is easy to see that for any $0<s<t$, on the event $\{M_{s,t}\le \delta\Lambda(t)\}$, it holds that $\cI (\varphi, \mathbb{Y}_t/\Lambda(t))=\cI (\varphi, \mathbb{Y}_{s,t}/\Lambda(t))$. Thus, for any $0<s<t$,
	\begin{align}\label{631}
\left |\E\left(\cI (\varphi, \mathbb{Y}_t/\Lambda(t))\right)-\E\left(\cI (\varphi, \mathbb{Y}_{s,t}/\Lambda(t))\right)\right|\le \P(M_{s,t}>\delta\Lambda(t)).
	\end{align}	
	Using the Markov property and  \eqref{314}, we have that
	\begin{align} \label{641}
		  &1-\E\left(\cI (\varphi, \mathbb{Y}_{s,t}/\Lambda(t))|\FF_{t-s}\right)
		 = 	 1-\prod\limits_{v \in \LL_{t-s}}\E\left(\cI (\varphi, \mathbb{Y}^v_{s,t}/\Lambda(t))|\FF_{t-s}\right)\nonumber\\
		 =
		&	1-\left[\E\left(\cI (\varphi, \mathbb{Y}'_s /\Lambda(t))\right)\right]^{Z_{t-s}}.
	\end{align}	
		It follows from Lemma \ref{lem:Ms} that, as $t\to\infty$,
	$$
	1-\E\left(\cI (\varphi, \mathbb{Y}'_s /\Lambda(t))\right)\sim C_s\Lambda(t)^{-\alpha}
	L(\Lambda(t)),
	$$
	where $C_s=\int_0^s \left(e^{\lambda (s-r)}-e^{-\beta (s-r)}\right)\int_{\R}\E \left(1-\varphi(y)^{Z_{r}}\right)v_\alpha(dy)\,\d r$. Thus
	$$
	-\log \E\left(\cI (\varphi, \mathbb{Y}'_s /\Lambda(t))\right)\sim 1-\E\left(\cI (\varphi, \mathbb{Y}'_s /\Lambda(t))\right)\sim C_s\Lambda(t)^{-\alpha}L(\Lambda(t))
	$$
	and consequently
	$$
		      -Z_{t-s}\log \E\left(\cI (\varphi, \mathbb{Y}'_s /\Lambda(t))\right)\sim e^{\lambda(t-s)}C_s\Lambda(t)^{-\alpha}L(\Lambda(t))\cdot W,
	$$
	which tends to 0 as $t\to \infty$ since $\lim\limits_{t \to \infty} \dfrac{\Lambda(t)}{h(t)}=\infty$.
	Combining this with \eqref{641}, we get that, as $t\to\infty$,
	\begin{align}\label{limit-M2}
		&\lim\limits_{t \to \infty} 	e^{-\lambda t}\Lambda(t)^{\alpha}L(\Lambda(t))^{-1}	\left(1-\E\left(\cI (\varphi, \mathbb{Y}_{s,t}/\Lambda(t))|\FF_{t-s}
		\right) \right) \nonumber\\
		 =&
		 \lim\limits_{t \to \infty} 	e^{-\lambda t}  \Lambda(t)^{\alpha}L(\Lambda(t))^{-1}
		 Z_{t-s}(-\log \E\left(\cI (\varphi, \mathbb{Y}'_s /\Lambda(t))\right))\nonumber  \\
		=&e^{-\lambda s} \int_0^s \left(e^{\lambda (s-r)}-e^{-\beta (s-r)}\right)\int_{\R}\E \left(1-\varphi(y)^{Z_{r}}\right)v_\alpha(\d y)\,\d r\cdot W.
	\end{align}	
	Moreover,
	\begin{align}\label{3}
		&e^{-\lambda t}\Lambda(t)^{\alpha}L(\Lambda(t))^{-1}	\left(1-\E\left(\cI (\varphi, \mathbb{Y}_{s,t}/\Lambda(t))|\FF_{t-s}
		\right) \right) \nonumber\\
		\le & 	e^{-\lambda t} Z_{t-s} \Lambda(t)^{\alpha}L(\Lambda(t))^{-1}  \left( 1-\E\left(\cI (\varphi, \mathbb{Y}'_{s}/\Lambda(t))\right)\right)\nonumber\\
		\overset{t}{\le}&2e^{-\lambda s} \int_0^s \left(e^{\lambda (s-r)}-e^{-\beta (s-r)}\right)\int_{\R}\E \left(1-\varphi(y)^{Z_{r}}\right)v_\alpha(\d y)\,\d r\cdot
			e^{-\lambda (t-s)}Z_{t-s}.
	\end{align}	
	In the first inequality we used the inequality	$1-x^n\leq n(1-x) , x \in(0,1)$.
Note that
	\begin{align}\label{W-{t-s}}
			e^{-\lambda (t-s)}Z_{t-s}
		\to W \quad \text{a.s.} , \quad \mbox{and }
			\E \left(e^{-\lambda (t-s)}Z_{t-s}\right) \to \E W=1,     \quad t \to \infty.
	\end{align}
	Combining  \eqref{limit-M2}, \eqref{3} and \eqref{W-{t-s}},  and using the dominated convergence theorem we get
	\begin{align}\label{000}
	&\lim\limits_{t \to \infty} 	e^{-\lambda t}\Lambda(t)^{\alpha}L(\Lambda(t))^{-1}	\left(1-\E\left(\cI (\varphi, \mathbb{Y}_{s,t}/\Lambda(t))
	\right) \right) \nonumber\\
	=&e^{-\lambda s} \int_0^s \left(e^{\lambda (s-r)}-e^{-\beta (s-r)}\right)\int_{\R}\E \left(1-\varphi(y)^{Z_{r}}\right)v_\alpha(\d y)\,\d r.
	\end{align}	
	Using  \eqref{631}, Lemma \ref{lem:m1} and \eqref{000},  we get that
	\begin{align*}
		&\lim\limits_{t \to \infty} 	e^{-\lambda t}\Lambda(t)^{\alpha}L(\Lambda(t))^{-1}	\left(1-\E\left(\cI (\varphi, \mathbb{Y}_{t}/\Lambda(t))
		\right) \right) \nonumber\\
		=& \lim_{s\to\infty}e^{-\lambda s} \int_0^s \left(e^{\lambda (s-r)}-e^{-\beta (s-r)}\right)\int_{\R}\E \left(1-\varphi(y)^{Z_{r}}\right)v_\alpha(\d y)\,\d r\\
		=& \int_0^\infty e^{-\lambda r}\int_{\R}\E \left(1-\varphi(y)^{Z_{r}}\right)v_\alpha(\d y)\,\d r=C(\varphi).
	\end{align*}	

We now prove \eqref{lim-Mt}. For any $\epsilon\in(0,1),$
define $\tilde{I}_{\varepsilon}, I_{\varepsilon} \in\mathcal{H}(\R)$ by
$$ I_{\varepsilon}(y):=\left\{\begin{array}{cc}
	1,&y<1, \\
	{\rm linear}, &1\leq y\leq 1+\varepsilon, \\
	0,&y>1+\varepsilon
\end{array}\right.
$$
and
$$\tilde{ I}_{\varepsilon}(y):=\left\{\begin{array}{cc}
	1,& y<1 -\varepsilon, \\
	{\rm linear}, &1-\varepsilon\leq y\leq 1, \\
	0,& y>1.
\end{array}\right.$$
It is clear that
 $ \tilde{I}_{\varepsilon}(y)\le  {\bf 1}_{(-\infty, 1]}(y)\leq  I_{\varepsilon}(y). $
Applying \eqref{lim-Yt} to $\tilde{I}_{\varepsilon}$ and $ I_{\varepsilon}, $ and by the continuity of $C(\varphi)$,  we can get  \eqref{lim-Yt} still holds for $\varphi(y)={\bf 1}_{(-\infty, 1]}(y)$, that is,
\begin{align*}
	\lim\limits_{t\to \infty} e^{-\lambda t}\Lambda(t)^{\alpha}L(\Lambda(t))^{-1} \P(M_t > \Lambda(t))&=C({\bf 1}_{(-\infty, 1]})\\
	&=\int_0^\infty e^{-\lambda r}\int_1^\infty \P(Z_r>0)v_\alpha(\d y) \d r=\frac{q_1}{\alpha}\vartheta.
\end{align*}
	The proof is now complete.
	
	$\hfill\qedsymbol$

	\noindent \textbf{Proof of Theorem \ref{them:large-deviation}}	
	Applying Lemma \ref{lem:one-big'} and Proposition \ref{prop:large-M}, we have that
	 \begin{align}\label{646}
		\lim\limits_{t\to \infty} e^{-\lambda t}\Lambda(t)^{\alpha}L(\Lambda(t))^{-1}\left( 1-\E\cI (\varphi, \mathbb{X}_t/\Lambda(t))\right)  =C(\varphi),\quad \varphi\in\mathcal{H}(\R).
	\end{align}
By the continuity of $C(\varphi)$,  we can get  \eqref{646} still holds for $\varphi(y)={\bf 1}_{(-\infty, 1]}(y)$, that is,
\begin{align}\label{111}
	\lim\limits_{t\to \infty} e^{-\lambda t}\Lambda(t)^{\alpha}L(\Lambda(t))^{-1} \P(R_t > \Lambda(t))&=C({\bf 1}_{(-\infty, 1]})=\vartheta^*.
\end{align}

\hfill$\Box$
	
	\noindent{\bf Proof of Corollary \ref{cor:large-deviation}:}
By \eqref{111}, we have  for $x>1$,
	\begin{align}
		\lim\limits_{t\to \infty}  \P(R_t>x\Lambda(t)|R_t > \Lambda(t))=x^{-\alpha}=P(R^*>x).
	\end{align}
For any $g\in C_c^+(\overline{\R}_0)$ and $\theta>0$, applying Theorem \ref{them:large-deviation}  with $\varphi(y)=e^{-\theta g(y)}{\bf 1}_{(-\infty, 1]}(y)$  and $\varphi(y)=e^{-\theta g(y)}$, we have that
	\begin{align}\label{112}
		&\lim\limits_{t\to \infty} e^{-\lambda t}\Lambda(t)^{\alpha}L(\Lambda(t))^{-1}\E\Big(e^{-\theta \frac{\mathbb{X}_t}{\Lambda(t)}(g)},R_t> \Lambda(t)\Big)\nonumber\\
		=&\lim\limits_{t\to \infty} e^{-\lambda t}\Lambda(t)^{\alpha}L(\Lambda(t))^{-1}\left[1-\E\Big(e^{-\theta \frac{\mathbb{X}_t}{\Lambda(t)}(g)},R_t\le \Lambda(t)\Big)-\Big(1-\E\Big(e^{-\theta \frac{\mathbb{X}_t}{\Lambda(t)}(g)}\Big)\Big)\right]\nonumber\\
		=&\int_0^\infty e^{-\lambda r}\int_{1}^\infty \E[e^{-\theta g(y)Z_r};Z_r>0] v_\alpha(\d y) \d r.
	\end{align}
	By \eqref{111} and \eqref{112}, we have
	\begin{align*}
		\lim_{t\to\infty}\E\Big(e^{-\theta \frac{\mathbb{X}_t}{\Lambda(t)}(g)}|R_t> \Lambda(t)\Big)&=\frac{\alpha }{q_1  \vartheta     }\int_0^\infty e^{-\lambda r}\int_{1}^\infty \E[e^{-\theta g(y)Z_r};Z_r>0] v_\alpha(\d y) \d r\\
	&=E\left(e^{-\theta T\cdot g(R^*)}\right).
	\end{align*}
	The proof is now complete. \hfill$\Box$

\section{Lower deviation of $\X_t$ and $\Y_t$}
	
	We first give some results about the continuous time branching processes $\{Z_t:t\ge0\}$.
	Recall that the extinction probability $q\in [0, 1)$.
	For any $s\in[0,1]$ and $t\ge0$, define
	$$
	F(s,t):=\E\left(s^{Z_t}\right).
	$$
	Recall the constant $\rho$ defined in \eqref{def:rho}.
	 It is well known  (see,
	 \cite[Section III.8]{Athreya-Ney}, for instance) that
	\begin{align}\label{limit-F}
		\lim_{t\to\infty} e^{\rho t}[F(s,t)-q]
		=:A(s)
	\end{align}
	exists for $0\le s<1$. Moreover, the convergence is uniform in $s\in [0,a]$ for any $0<a<1$. The function $A(s)$ is the unique solution of
	\begin{align}\label{622}
		A(F(s,t))=e^{-\rho t}A(s)
	\end{align}
	with $A(q)=0, A'(q)=1$.
	Since $A(s)$ is the limit of power series, it is itself  a power series
	$$A(s)=\sum_{j=0}^\infty a_j s^j.$$
	It is clear that $a_0\le0$ and $a_j\ge 0, j\ge1$. For any $0<c<1-q$ and $s\in[0,1]$,
	\begin{align}
		A(cs+q)=A(cs+q)-A(q)=\sum_{j=1}^\infty a_j[(cs+q)^j-q^j]=\sum_{k=1}^\infty \Big[\sum_{j=k}^\infty a_j
\binom{j}{k}		
q^{j-k}\Big]c^ks^k.
	\end{align}
Thus $\frac{A(cs+q)}{A(c+q)}$,  $s\in[0,1]$, is a probability generating function.

	\begin{lemma} \label{lem:As}
		For any $s\in[0,1)$,
		$$A(s)=s-q+\int_0^\infty \beta e^{\rho t} V(F(s,r))\,\d r,$$
		where \begin{align}\label{est:V}
			0\le V(s):=f(s)-f'(q)s-q(1-f'(q))\le f''(s\vee q)(s-q)^2.
		\end{align}
	\end{lemma}
	
	\noindent{\bf Proof:}
	By the Markov property and the branching property, we have that
	\begin{align*}
		F(s,t)=&se^{-\beta t}+\int_0^t \beta e^{-\beta r} f(F(s,t-r))\d r\\
		=&se^{-\rho t}+\int_0^t \beta e^{-\rho r} f(F(s,t-r))\d r-\int_0^te^{-\rho r} \beta f'(q)F(s,t-r)\d r \\
		=&se^{-\rho t}+\int_0^t \beta e^{-\rho r} V(F(s,t-r))\d r+q\rho\int_0^t e^{-\rho r}\d r\\
		=&e^{-\rho t}\left(s-q+\int_0^t  \beta e^{\rho r}V(F(s,r))\,\d r \right) +q.
	\end{align*}
By \eqref{limit-F}, we have
	\begin{align*}
		A(s)&=\lim_{t\to\infty}e^{\rho t}(F(s,t)-q)=s-q+\lim_{t\to\infty}\int_0^t  \beta e^{\rho r}V(F(s,r))\,\d r\\
		&=s-q+\int_0^\infty \beta e^{\rho r}V(F(s,r))\,\d r.
	\end{align*}
Since $q=f(q)$,  we have
	$
	V(s)=f(s)-f(q)-f'(q)(s-q).
	$
	Now \eqref{est:V} follows immediately from Taylor's formula.
	\hfill$\Box$.

		Recall that $\phi(\theta)=\E(e^{-\theta W})$.
	\begin{lemma}
		For any $\theta>0$,
		\begin{align*}
		A\left(\phi(\theta )\right)=	\int_{-\infty}^\infty\beta e^{\rho s} V \left[\phi\left( \theta e^{\lambda s}\right) \right]\,ds.
		\end{align*}
		
	\end{lemma}
	
	\noindent	{\bf Proof:}
	By the branching property and the Markov property, we have
	$$Z_{t+s}=\sum_{u\in\LL_{s}}Z^{u}_{t+s}.$$
	Given $Z_s$, $\{Z^{u}_{t+s},u\in\LL_{s}\}$ are  i.i.d. with the same law as $Z_t$. It follows that
	$$W=\lim_{t\to\infty}e^{-\lambda s}\sum_{u\in\LL_{s}}e^{-\lambda t}Z^{u}_{t+s}=e^{-\lambda s}\sum_{u\in\LL_{s}}W^{u},$$
	where  $\{W^{u},u\in\LL_{s}\}$ are i.i.d. with the same law as $W$.
	Thus we have
	\begin{align}\label{623}
		\phi(\theta e^{\lambda s})=\E\left( \phi(\theta)^{Z_s}\right) =F( \phi(\theta),s).
	\end{align}
	Hence by Lemma \ref{lem:As} and \eqref{623}, we have for any $x>0$
	\begin{align*}
		& \int_{-x}^\infty\beta e^{\rho s} V \left[\phi\left( \theta e^{\lambda s}\right) \right]\,ds
		= e^{-\rho x}\int_{0}^\infty\beta e^{\rho s} V \left[(\phi(\theta e^{\lambda s}e^{-\lambda x} \right]\,ds\\
		=&e^{-\rho x}\int_{0}^\infty\beta e^{\rho s} V \left[F(\phi(\theta e^{-\lambda x}),s) \right]\,ds=e^{-\rho x}[A\left(\phi(\theta e^{-\lambda x}) \right)- \phi(\theta e^{-\lambda x})+q]\\
		=&A(F(\phi(\theta e^{-\lambda x}),x))-e^{-\rho x}(\phi(\theta e^{-\lambda x})-q)\\
		=&A(\phi(\theta))-e^{-\rho x}(\phi(\theta e^{-\lambda x})-q),
	\end{align*}
where the fourth equality follows from \eqref{622}.
	Letting $x\to\infty$, the desired result follows immediately.
	\hfill$\Box$

	\begin{lemma}\label{lem:lambda-rho}
		It holds that  $\lambda\ge \rho$.
		Moreover, $\lambda=\rho$ if and only if $p_k=0$ for all $k\ge3$.
	\end{lemma}
	\noindent{\bf Proof:} Note that $\lambda-\rho=\beta(f'(1)+f'(q)-2)$.
	
We first consider the case when $p_0=0$. Then $q=0$ and $f'(q)=f'(0)=p_1.$
	Hence we have
	$$f'(1)+f'(0)-2= \left( \sum_{k=1}^\infty kp_k+p_1-2\sum_{k=1}^\infty p_k\right)=\sum_{k=2}^\infty(k-2)p_k\ge0. $$
Furthermore,  the left hand side is equal to 0  if and only if $p_k=f^{(k)}(0)=0$, $k\ge3$.
	
	If $p_0>0$, then $q>0$. Define
	$$\hat{f}(s):=\frac{f((1-q)s+q)-q}{1-q},\quad \mbox{for } 0\le s\le 1.$$
	It is well known (see \cite[Chapter 1, Section 12]{Athreya-Ney}) that $\hat{f}$ is a probability generating function with $\hat{f}(0)=0$.
	Note that $\hat{f}'(1)=f'(1)$ and $\hat{f}'(0)=f'(q)$.
	Applying the previous paragraph to $\hat{f}$, we get that
	$$f'(1)+f'(q)-2=\hat{f}'(1)+\hat{f}'(0)-2\ge0.$$
	Moreover,  the left hand side is equal to 0
	if and only if $\hat{f}^{(k)}(0)=(1-q)^{k-1}f^{(k)}(q)=0, k\ge3.$
	It is easy to see that $f^{(3)}(q)=\sum_{k=3}^\infty k(k-1)(k-2)p_kq^{k-3}=0$ if and only if $p_k=0$, $k\ge3$.
	
The proof is now complete.

	\hfill$\Box$

	Note that  the skeleton $\{Z_n,n=0,1,\cdots\}$ is a Galton-Watson process with offspring  generating function $\tilde{f}(s):=F(s,1)$.  It is clear that  $\tilde{f}'(1)=
	\E(Z_1)=e^{\lambda}>1$ and $\P(Z_1=1)\ge\P(\tau_o>1)>0$ .
	By \cite[Theorem 4]{FW}, we have that,	for any nonnegative sequence $a_n$ with $a_n\to\infty$ and $a_n=o(e^{\lambda n})$,  there exists a constant $c>0$ such that
	\begin{align}\label{351}
		\P(0<Z_n< a_n)
		\le c\P(0<W<e^{-\lambda n}a_n),\qquad n\ge1.
	\end{align}
	Furthermore,
	it was proven  in \cite{Dubuc} that there exist $c_1, c_2>0$ such that  for any $x\in(0,1)$,
	\begin{align}\label{352}
		c_1 x^{\rho/\lambda}<\P(0<W<x)<c_2 x^{\rho/\lambda}.
	\end{align}
	By \eqref{351} and \eqref{352},   there exists $C>0$ such that
	\begin{align}\label{est-Zn}
		\P(0<Z_n< a_n)
		\le  C a_n^{\rho/\lambda}e^{-\rho n}, \quad n\ge1.
	\end{align}

	\subsection{Lower deviation of $\X_t$}\label{subsec3.1}

We have proved in \cite[Proposition 2.1]{YRR} that 	$\mathbb{Y}_t/h(t)$ converges weakly to
$\mathcal{N}_\infty$, that is for any $\varphi=e^{-g}$ with $g\in C_c^+(\overline{\R}_0)$,
\begin{align}\label{laplace-N3}
	\lim_{t\to\infty}\E(\cI(\varphi, \mathbb{Y}_t/h(t)))=\E((\cI(\varphi, \mathcal{N}_\infty))
	&=\E\left(\exp\left\{-C(\varphi)W\right\}\right).
\end{align}
In fact, by examining the proofs in \cite{YRR},  \eqref{laplace-N3} is valid  for any $\varphi\in \cH(\R)$.
By Lemma \ref{lem:one-big'} we have for any $\varphi\in \cH(\R)$,
	\begin{align}\label{laplace-N2}
		\lim_{t\to\infty}\E(\cI(\varphi, \mathbb{X}_t/h(t)))=\E((\cI(\varphi, \mathcal{N}_\infty))
		&=\E\left(\exp\left\{-C(\varphi)W\right\}\right).
	\end{align}
For any function $g$ and $x\in\R$,
we define the function $\mathfrak{m}_x g(\cdot)$ by the relation  $\mathfrak{m}_x g(y)=g(y/x)$.
By the definition of $C(\varphi)$ in \eqref{def:C}, it is easy to see that
\begin{align}
	C(\mathfrak{m}_x\varphi)=x^{-\alpha}C(\varphi),\quad\varphi\in\mathcal{B}_1(\R).
\end{align}

\begin{lemma}\label{lem:11}
If $a(t), b(t)$ are positive functions with $\frac{a(t)}{h(t)}\to a>0$ and $\frac{b(t)}{a(t)}\to0$, then
for any $\varphi\in \cH(\R)$,
	\begin{align}
		\lim_{t\to\infty}\E_{b(t)}(\cI(\varphi, \mathbb{X}_t/a(t)))=\E(\exp\{-a^{-\alpha} C(\varphi)W\}).
		\end{align}	
	\end{lemma}

\noindent{\bf Proof:}  Note that $$\E_{b(t)}(\cI(\varphi, \mathbb{X}_t/a(t)))=\E\left(\prod_{u\in\mathcal{L}_t}\varphi\left(\frac{\xi_t^u+b(t)}{a(t)}\right)\right).$$

Since $\varphi\in\mathcal{H}(\R)$, we have
$\varphi\equiv1$ on $[-\delta, \delta]$ for some $\delta>0$.
Moreover for any $\epsilon>0$, there exists $\eta=\eta(\epsilon)>0$ such that $|\varphi(x)-\varphi(y)|<\epsilon$ whenever
$|x-y|\le \eta$.
It follows from the assumption that for any $\epsilon'>0,$ there exists $t'$ such that for all $t>t'$,
$(1-\epsilon')a(t)\le a h(t)\le (1+\epsilon')a(t)$
and $|b(t)|\le \epsilon' a(t).$

We now fix an arbitrary $\epsilon>0$ and $0<\epsilon'<\frac{\delta}{\delta+2}\wedge \eta$.
Note that if $\frac{|y|}{ah(t)}\le \delta/2$,
then for $t>t'$,
$$
\left|\frac{y+b(t)}{a(t)}\right|\le \frac{(1+\epsilon')|y| }{ah(t)} +\epsilon' \le\left(1+\frac{\delta}{\delta+2}\right)\frac{\delta}{2}+\frac{\delta}{\delta+2}= \delta,$$
 and that if
  $\frac{|y|}{ah(t)}\le  \frac{\eta}{\epsilon'}-1$, then for $t>t'$,
 $$\left|\frac{y+b(t)}{a(t)}-\frac{y}{a h(t)}\right|\le \epsilon'+\frac{|y|}{a h(t)}\epsilon'\le \eta.$$
Thus for any $t>t'$,
 $$\left|\varphi\Big(\frac{y+b(t)}{a(t)}\Big)-\varphi\Big( \frac{|y|}{ah(t)}\Big)\right|\le
 \epsilon {\bf 1}_{\{\frac{\delta}{2}<\frac{|y|}{ah(t)}\le \frac{\eta}{\epsilon'}-1\}}+ {\bf 1}_{\{ \frac{|y|}{ah(t)}> \frac{\eta}{\epsilon'}-1\}}.
 $$
 Hence we have that
 \begin{align}\label{651}
& \left|\E\left(\prod_{u\in\mathcal{L}_t}\varphi\left(\frac{\xi_t^u+b(t)}{a(t)}\right)\right)-\E\left(\prod_{u\in\mathcal{L}_t}\varphi\left(\frac{\xi_t^u}{a h(t)}\right)\right)\right|\nonumber\\
\le& \E\left[1\wedge \sum_{u\in\mathcal{L}_t}\left|\varphi\left(\frac{\xi_t^u+b(t)}{a(t)}\right)-\varphi\left(\frac{\xi_t^u}{a h(t)}\right)\right|\right]\nonumber\\
 \le& \E\left[1\wedge \left(\epsilon \sum_{u\in\mathcal{L}_t} {\bf 1}_{\{\frac{|\xi_t^u|}{h(t)}>a\delta /2\}}\right)\right]+\E\left[1\wedge \sum_{u\in\mathcal{L}_t} {\bf 1}_{\{\frac{|\xi_t^u|}{h(t)}>a(\eta/\epsilon'-1)\}}\right].
 \end{align}
Since $\mathbb{X}_t/h(t)$ converges weakly to $\mathcal{N}_\infty$,
we have, as $t\to\infty$,
 $$\epsilon \sum_{u\in\mathcal{L}_t} {\bf 1}_{\{\frac{|\xi_t^u|}{h(t)}>a\delta /2\}}\overset{d}{\to }
\epsilon\cdot  \mathcal{N}_\infty(\{y\in\R; |y|>a\delta/2\})
 $$
 and
$$\sum_{u\in\mathcal{L}_t} {\bf 1}_{\{\frac{|\xi_t^u|}{h(t)}>a(\eta/\epsilon-1)\}}\overset {d}{\to }
\mathcal{N}_\infty(\{y\in\R; |y|>a(\eta/\epsilon'-1)\}.
$$
Thus letting
$t\to\infty$ first and then $\epsilon'\to 0$,  and finally $\epsilon\to0$ in \eqref{651}, applying the dominated convergence theorem, we get
\begin{align*}
&\lim_{t\to\infty}\E\left(\prod_{u\in\mathcal{L}_t}\varphi\left(\frac{\xi_t^u+b(t)}{a(t)}\right)\right)=\lim_{t\to\infty}\E\left(\prod_{u\in\mathcal{L}_t}\varphi\left(\frac{\xi_t^u}{a h(t)}\right)\right)\\
&=\E\left(\exp\left\{-C(\mathfrak{m}_a\varphi)W\right\}\right)=\E(\exp\{-a^{-\alpha} C(\varphi)W\}).
\end{align*}

\hfill$\Box$

	\begin{lemma} \label{lem:limit-xi}
	If   $a:[0,\infty)\to(0,\infty)$ is a
	non-decreasing positive function
	with $$\sum_{n=1}^\infty n a(n)^{-\alpha}L(a(n))<\infty,$$ then
	$$\lim_{t\to\infty}\frac{\xi_t}{a(t)}=0.$$
\end{lemma}
\noindent{\bf Proof:}
For any positive integer $n$,
$$\sup_{t\in[n,n+1]}|\xi_t|\le |\xi_n|+\sup_{t\in[n,n+1]}|\xi_t-\xi_n|.$$
Since $\xi$ is a L\'evy process,
$\{Y_n:=\sup_{t\in[n,n+1]}|\xi_t-\xi_n|, n\ge1 \}$ are  i.i.d.
By \eqref{sup-xi}, and \eqref{inf-xi},  we have that as $x\to\infty$,
$$\mP(Y_1>x)\le \mP\left(\sup_{t\in[0,1]}\xi_t>x\right)+\mP\left(\sup_{t\in[0,1]}(-\xi_t)>x\right)\sim \mP(|\xi_1|>x).$$
By \eqref{tail-prob}, we have that for any $c>0$ and $n$ large enough,
\begin{align*}
	\mP\left(\sup_{t\in[n,n+1]}|\xi_t|>ca(n)\right)\le & \mP(|\xi_n|>ca(n)/2)+\mP(Y_n>ca(n)/2)\\
	\le &c_0(c/2)^{-\alpha} (n+1)a(n)^{-\alpha}L(ca(n)/2)\\
	\sim &c_0(c/2)^{-\alpha} na(n)^{-\alpha}L(a(n)), \quad n\to\infty,
\end{align*}
which implies that
$$\sum_n\mP\left(\sup_{t\in[n,n+1]}|\xi_t|>ca(n)\right)<\infty.$$
Thus by the Borel-Cantelli lemma, we have
$$\lim_{n\to\infty}\frac{\sup_{t\in[n,n+1]}|\xi_t|}{a(n)}=0,\quad a.s.$$
Since $a(t)$ is non-decreasing, we have, for $t\in[n,n+1]$,
$$\frac{|\xi_t|}{a(t)}\le \frac{\sup_{t\in[n,n+1]}|\xi_t|}{a(n)}\to0.$$
The proof is now complete.

\hfill$\Box$

For any $\varphi\in\mathcal{B}_1(\R)$, we define
\begin{align}\label{d:Uvarphi}
U_\varphi(t,x):=\E_x(\mathcal{I}(\varphi,\X_t))=\E_{x}\left( \prod_{u\in\mathcal{L}_t}\varphi(\xi_t^u)\right),\quad t\ge0,
x\in \R.
\end{align}
By the Markov property and the branching property, we have
$$U_\varphi(t,x)=
	e^{-\beta t}\mE_{x}(\varphi(\xi_t))+\mE_{x}\int_0^t \beta e^{-\beta s} f(U_\varphi(t-s, \xi_s))ds,$$ which implies that
\begin{align*}
	U_\varphi(t,x)
	&=	e^{\beta(f'(q)-1) t}\mE_{x}(\varphi(\xi_t))+\mE_{x}\int_0^t \beta e^{\beta(f'(q)-1) s} f(U_\varphi(t-s, \xi_s))\,ds\\
	&\quad -\mE_x\int_0^te^{\beta(f'(q)-1) s} \beta f'(q)U_\varphi(t-s,\xi_s)\,ds\\
	&=e^{-\rho t}\mE_x(\varphi(\xi_t))+\mE_x\int_0^t \beta e^{-\rho s} V(U_\varphi(t-s, \xi_s))\,ds+q(1-e^{-\rho t}),
\end{align*}
where $\rho$ is defined in \eqref{def:rho} and $V$ is defined in \eqref{est:V}. Thus
\begin{align}\label{eq:u}
	U_\varphi(t,x)-q=e^{-\rho t}(\mE(\varphi(\xi_t+x))-q)+\mE\int_0^t \beta e^{-\rho s} V(U_\varphi(t-s,x+ \xi_s))\,ds.
\end{align}

For any $t>0$ and $x\in\R$, we define
\begin{align}\label{d:u}
u(t,x):=U_{{\bf 1}_{(-\infty,0]}}(t,x)=\P_x(R_t\le 0)=\P(R_t\le -x).
\end{align}

	\begin{lemma}\label{lem:est-u}
		Let $c>0$. For any $\epsilon\in(0,1)$, there exist $C=C(\epsilon)>0$ and $t_\epsilon>0$ such that for any $r>t_\epsilon, s>t_\epsilon$,  $l>0$,
		$$\mE[(u(r+s,-c\,h(r)+\xi_{l})-q)^2]\le c(\epsilon)s^{6\rho/\lambda}\Big[ e^{-2\rho (1-\epsilon)s}+ (l+s) e^{-\lambda r} e^{-(\lambda\wedge (2\rho))(1-\epsilon)s}\Big].$$
		
	\end{lemma}
	\noindent{\bf Proof:}
	Define $\tilde{u}(t,x):=\P(-\infty<R_t\le -x)=u(t,x)-\P(Z_t=0)$. Then for any $s>0$,
\begin{align}\label{sum}
		&\mE[(u(r+s,-c\,h(r)+\xi_{l})-q)^2]
		\le 2\mE[(\tilde{u}(r+s,-c\,h(r)+\xi_{l}))^2]+2(q-\P(Z_{r+s}=0))^2\nonumber\\
		&\le 2\mE[(\tilde{u}(r+s,-c\,h(r)+\xi_{l})^2]+2(q-\P(Z_{s}=0))^2.
\end{align}
	By \eqref{limit-F} with $s=0$, we have
	\begin{align}\label{655}
		\lim_{s\to\infty} e^{\rho s}
		(q-\P(Z_{s}=0))=-A(0)\in[0,\infty).
	\end{align}
	
		For $t>0$, define
	$$T_t:=\inf\{s>0:Z_{s}>t^3\}.$$
	Here we use the convention that $\inf \emptyset=\infty$.
Let $K>3$ be an integer.
For any $r>0$, $s>1$ and $x\in\R$, we have
	\begin{align}\label{613}
		&\tilde{u}(s+r,-c\,h(r)+x)=\P(-\infty<R_{s+r}\le c\,h(r)-x)\nonumber\\
		&\le \P\left( T_s>\frac{K-1}{K} s,Z_{(K-1)s/K}>0\right)\nonumber\\
		\nonumber\\&
		\quad+\sum_{k=1}^{K-1}\P\left(\frac{(k-1)s}{K}< T_s\le \frac{ks}{K}, R_{s+r}\le c\,h(r)-x\right).
	\end{align}
	By \eqref{est-Zn} with $a_n=(2K)^3 n^3$ , there exists a constant $C_K>0$ such that
	\begin{align*}
		\P(0<Z_n\le (2K)^3 n^3)\le C_K (2Kn)^{3\rho/\lambda} e^{-\rho n}, \quad n\ge1.
	\end{align*}
	It follows that, for $s>K$ and $k=1,2,\cdots, K-1$,
	\begin{align}\label{612}
		\P\left( T_s>\frac{k}{K} s,Z_{ks/K}>0\right)&\le \P\left( 0<Z_{\lfloor\frac{k}{K}s\rfloor}\le s^3\right) \le \P\left( 0<Z_{\lfloor\frac{k}{K}s\rfloor}\le (2K)^3 \lfloor\frac{k}{K}s\rfloor^3\right)\nonumber\\
		&\le  C_K \left(2Ks\right)^{3\rho/\lambda}e^{-\rho\lfloor\frac{k}{K}s\rfloor}\le  C_K (2K)^{3\rho/\lambda} e^{\rho} \cdot s^{3\rho/\lambda} e^{-\rho \frac{k}{K}s}.
	\end{align}

Choose $b>1$ such that $e^{-\lambda b/\alpha}<1/c$.	
Note that $\lim_{t\to\infty}\frac{h(t)}{h(t+b)}=e^{-\lambda b/\alpha}<1/c$.
By \eqref{Rt/ht},  we have $\lim_{t\to\infty}\P(R_t\le h(t))=\E(e^{-\vartheta^* W})<1$. Thus there exist $c_1\in (0, 1)$, $c_2\in (0, 1/c)$  and $t_1>0$
such that  for all $t>t_1$, and $t'>b$,
	$$\P(R_t\le h(t))\le c_1,\quad \mbox{and } h(t)\le c_2h(t+b)\le c_2h(t+t') .$$
	
	Note that $$R_{s+r}=\max_{u\in\mathcal{L}_{ks/K}}(\xi_{ks/K}^u+R^u_{s+r}),$$
	where,  given $\FF_{sk/K}^{\mathbb{T}}$, $\{R^u_{s+r}, u\in\mathcal{L}_{ks/K}\}$  are i.i.d. with the same law as $R_{s+r-ks/K}$. Moreover,  given $\FF_{ks/K}^{\mathbb{T}}$,  $\{\xi_{ks/K}^u,u\in\mathcal{L}_{ks/K}\}$ have the same law as $\xi_{ks/K}$, and $\{\xi_{ks/K}^u,u\in\mathcal{L}_{ks/K}\}$  and $\{R^u_{s+r}, u\in\mathcal{L}_{ks/K}\}$  are independent.
Applying \cite[Lemma 5.1]{GH} (in the first inequality below) we get that, for $k=1,2,\cdots, K-1$,  $r>t_1$ and $s>bK$,
	\begin{align*}
		&\P(R_{s+r}\le c\,h(r)-x|\mathcal{F}^{\mathbb{T}}_{ks/K})\le \P(\xi_{ks/K}+\max_{u\in\mathcal{L}_{ks/K}}R^u_{s+r}\le c\,h(r)-x|\mathcal{F}^{\mathbb{T}}_{ks/K})\nonumber\\
		\le &\mP(x+\xi_{ks/K}\le c\,h(r)-h(s+r-ks/K))
		+[\P(R_{s+r-{ks/K}}\le h(s+r-{ks/K}))]^{Z_{ks/K}}\nonumber\\
		\le& \mP(x+\xi_{ks/K}\le -(1-cc_2)h(s+r-ks/K))
		+c_1^{Z_{ks/K}}.
	\end{align*}
	Thus we have that for $k=1,2,\cdots,K-1$, $r>t_1$ and $s>bK$,
	\begin{align}\label{671}
		&\P\left(\frac{(k-1)s}{K}< T_s\le \frac{ks}{K}, R_{s+r}\le c\,h(r)-x\right)\nonumber\\
		\le&  \mP(x+\xi_{ks/K}\le -(1-cc_2)h(s+r-ks/K))\P\left(T_s>\frac{(k-1)s}{K}, Z_{(k-1)s/K}>0\right)\nonumber\\
		&+
		\E
		\left( c_1^{Z_{ks/K}}; \frac{(k-1)s}{K}<T_s\le \frac{ks}{K}\right).
	\end{align}
	Note that
	\begin{align}\label{621}
		&
		\E\left( c_1^{Z_{ks/K}}; \frac{(k-1)s}{K}< T_s\le \frac{ks}{K}\right)\nonumber\\
		\le& c_1^{s^2}+\P\left(  \frac{(k-1)s}{K}<T_s\le \frac{ks}{K}, Z_{ks/K}<s^2\right)\nonumber\\
		\le& c_1^{s^2}+\P(\exists t\ge0: Z_t\le \lfloor s^2\rfloor|Z_0=\lfloor s^3\rfloor)\nonumber\\
		\le& c_1^{s^2}+{{\lfloor s^3\rfloor}\choose{\lfloor s^2\rfloor}}q^{\lfloor s^3\rfloor-\lfloor s^2\rfloor}\le c_1^{s^2}+(s^3+1)^{s^2+1}q^{s^3-s^2-2}=o(e^{-\rho s}).
	\end{align}
		Combining \eqref{613}, \eqref{612}, \eqref{671} and
 \eqref{621}, we get that there exist $C_k'$ and $t_2>0$ such that for any $r>t_2, s>t_2$ and $x\in\R$,
	\begin{align}\label{626}
		&	[\P(-\infty<R_{s+r}\le h(r)-x)]^2 \le C'_K\Big[s^{6\rho/\lambda} e^{-2\rho \frac{K-1 }{K}s}\nonumber\\
		&\quad +\sum_{k=1}^{K-1} s^{6\rho/\lambda} e^{-2\rho(k-1)s/K}\mP(x+\xi_{ks/K}\le -(1-cc_2)h(s+r-ks/K))\Big].
	\end{align}
	By the Markov property of $\xi$ and \eqref{tail-prob2}, we get that there exists $t_3>t_2$ such that for $r>t_3$,
	\begin{align*}
		&\mE[\mP(x+\xi_{ks/K}\le -(1-cc_2)h(s+r-ks/K))|_{x=\xi_{l}}]\\
		=&\mP(\xi_{l+ks/K}\le -(1-cc_2)h(s+r-ks/K))\\
		\le & 2c_0(1-cc_2)^{-\alpha-1} (l+s)e^{-\lambda(r+s-ks/K)}.
	\end{align*}
	It follows from \eqref{626} that for any $s>t_3$, $r>t_3$ and $l>0$,
	\begin{align}
		&\mE[	\tilde{u}(s+r,-h(r)+\xi_{l})]^2\nonumber\\
		\le &C''_K\Big[s^{6\rho/\lambda} e^{-2\rho \frac{K-1 }{K}s}
		+s^{6\rho/\lambda} (l+s)\sum_{k=1}^{K-1}  e^{-2\rho(k-1)s/K}e^{-\lambda(r+s-ks/K)}\Big]\nonumber\\
		\le &C''_K\Big[s^{6\rho/\lambda} e^{-2\rho \frac{K-1 }{K}s}+
		Ks^{6\rho/\lambda} (l+s) e^{-\lambda r} e^{-(\lambda\wedge (2\rho))\frac{K-1}{K}s}\Big],
	\end{align}
	where
	$C''_K=C'_k(1+2c_0(1-cc_2)^{-\alpha-1})$
	and  in the last inequality we
	used
	$$2\rho(k-1)+\lambda(K-k)\ge (\lambda\wedge 2\rho)(K-1), \quad k=1,\cdots,K.$$
	
	For any $\epsilon>0$, we choose $K$ such that $1/K<\epsilon$,
	we get that there exists $c(\epsilon)>0$ such that
	$$
	2\mE[(\tilde{u}(r+s,-c\,h(r)+\xi_{l})^2]
	\le c(\epsilon)s^{6\rho/\lambda}\Big[ e^{-2\rho (1-\epsilon)s}+ (l+s) e^{-\lambda r} e^{-(\lambda\wedge (2\rho))(1-\epsilon)s}\Big].
	$$
	Combining this with \eqref{sum} and  \eqref{655}, we immediately get the desired result. The proof is now complete.

	\hfill$\Box$

\noindent	{\bf Proof of Theorem \ref{them:lower deviation}:}
Since $\varphi\in\cH_0(\R)$,
there exists $c>0$ such that $\varphi(y)=0$ for all $y>c$.
For any $t>0$, define $r=r(t)$ by $h(r(t))=\Lambda(t)$, that is,
\begin{equation}\label{def-r}
e^{-\lambda r(t)}=\Lambda(t)^{-\alpha}L(\Lambda(t)).
\end{equation}

	Since $\lim_{t\to\infty}\frac{\Lambda(t)}{h(t)}=0$ and $\Lambda(t)\to\infty$, we have $t-r(t)\to\infty$ and $r(t)\to\infty$ as $t\to\infty$.
	Thus for any $T>0$, there exists $t_0>0$ such that $t-r(t)>T$ and $r(t)>T$ for all $t\ge t_0$.
	In the remainder of this proof, we assume $t\ge t_0$.
	Note that
	\begin{align}\label{648}
		&(1-q)\E^*\left(\prod_{u\in\mathcal{L}_t}\varphi(\xi^u_t/\Lambda(t))\right)\nonumber\\
		=&\E\left(\prod_{u\in\mathcal{L}_t}\varphi(\xi^u_t/\Lambda(t))\right)-q+q-\E\left(\prod_{u\in\mathcal{L}_t}\varphi(\xi^u_t/\Lambda(t));\mathcal{S}^c\right).
		\end{align}
By \eqref{limit-F} with $s=0$, we have
	\begin{align}\label{647}
		0\le q-\E\left(\prod_{u\in\mathcal{L}_t}\varphi(\xi^u_t/\Lambda(t));\mathcal{S}^c\right)\le q-\P(Z_t=0)\sim -A(0) e^{-\rho t}.
	\end{align}
	By \eqref{eq:u},  we have that
	\begin{align}\label{616}
		\E\left(\prod_{u\in\mathcal{L}_t}\varphi(\xi^u_t/\Lambda(t))\right)-q
		&=e^{-\rho t}(\mE(\varphi(\xi_t/ \Lambda(t)))-q)+J(t),
	\end{align}
	where
	\begin{align}\label{exp-J}
		J(t)=&\left( \int_0^{t-r(t)-T}+\int_{t-r(t)-T}^{t-r(t)+T}+\int_{t-r(t)+T}^t\right)  \beta e^{-\rho s} \mE [V(U_{\varphi_t}(t-s, \xi_s))]\d s\nonumber\\
		=:&J_1(t, T)+J_2(t, T)+J_3(t, T),
	\end{align}
 $\varphi_t(y)=\varphi(y/\Lambda(t))$ and
$U_{\varphi}$ is defined in \eqref{d:Uvarphi}.
	It is easy to see that
	\begin{align}\label{672}
	e^{-\rho t}(\mE(\varphi(\xi_t/ \Lambda(t)))-q)=o(e^{-\rho(t-r(t))}).
	\end{align}
	We now deal with the three components of $J(t)$ separately.
We will show that
$$
 \lim_{T\to\infty}\limsup_{t\to\infty}e^{\rho(t-r(t))}J_1(t, T)=\lim_{T\to\infty}\limsup_{t\to\infty}e^{\rho(t-r(t))}J_3(t, T)=0.
$$

	(1) Since $0\le V(s)\le f(s)\le 1$, we have
	\begin{align}\label{large-int}
	e^{\rho(t-r(t))}J_3(t, T)\le e^{\rho(t-r(t))}\int_{t-r(t)+T}^\infty \beta e^{-\rho s} \,ds=\frac{1}{1-f'(q)} e^{-\rho T}.
	\end{align}
Hence $\lim_{T\to\infty}\limsup_{t\to\infty}e^{\rho(t-r(t))}J_3(t, T)=0$.

	(2)
	Note that
	\begin{align}\label{645}
		e^{\rho(t-r(t))}J_2(t, T)=\int_{-T}^{T}\beta e^{\rho s}\mE V(U_{\varphi_t}(s+r(t)),\xi_{t-r(t)-s})\,ds.
	\end{align}
	Since $\Lambda(t)=h(r(t))$, we have
	\begin{align}
		U_{\varphi_t}(s+r(t),\xi_{t-r(t)-s})=\E_{\xi_{t-r(t)-s}}\left(\prod_{u\in\mathcal{L}_{s+r(t)}}\varphi(\xi_{s+r(t)}^u/h(r(t)))\right).
	\end{align}
	By Lemma \ref{lem:limit-xi} and the fact that $\Lambda$ is non-decreasing, we have that for any $s\in[-T,T]$,
	$$\lim_{t\to\infty}\frac{|\xi_{t-r(t)-s}|}{h(s+r(t))}=\lim_{t\to\infty}\frac{|\xi_{t-r(t)-s}|}{e^{\lambda s/\alpha}\Lambda(t)}\le \lim_{t\to\infty} e^{-\frac{\lambda}{\alpha}s}\frac{|\xi_{t-r(t)-s}|}{\Lambda(t-r(t)-s)}=0.$$
	Note that  $h(r(t))/h(s+r(t))\to e^{-\frac{\lambda}{\alpha}s}.$ By Lemma \ref{lem:11}, we have  that for $s\in[-T,T]$,
	\begin{align}
		\lim_{t\to\infty} U_{\varphi_t}(s+r(t),\xi_{t-r(t)-s})=\E(\exp\{-e^{\lambda s}C(\varphi)W\})=\phi(e^{\lambda s}C(\varphi)).
	\end{align}
	Since $V(s)$ is bounded and continuous, by \eqref{645} and the dominated convergence theorem we have
	\begin{align}\label{mid-int}
		\lim_{t\to\infty}e^{\rho(t-r(t))}
		J_2(t, T)
		=\int_{-T}^T \beta e^{\rho s} V \left[\phi(e^{\lambda s}C(\varphi))\right]\d s.
	\end{align}
	
	(3)
	We now deal with $J_1(t, T)$. Recall $u$ is defined in \eqref{d:u}.
	Note that
	\begin{align}\label{673}
		&e^{\rho(t-r(t))}J_1(t, T)=\int_{T}^{t-r(t)}\beta e^{\rho s}\mE V(U_{\varphi_t}(s+r(t),\xi_{t-r(t)-s}))\,ds\nonumber\\
		\le &\int_{T}^{t-r(t)}\beta e^{\rho s}\mE [V(u(s+r(t),-c\,
		h(r(t))
		+\xi_{t-r(t)-s}))]\d s,
	\end{align}
where in the last inequality we used the fact that
$$U_{\varphi_t}(s+r(t),x)\le \P_{x}(R_{s+r(t)}\le c \, h(r(t)) )=u(s+r(t), -c\, h(r(t))+x).$$
	By Taylor's formula and the fact that $f'(s)$ is increasing, we have
	\begin{align}\label{e:V}
	0\le V(s)=f(s)-f(q)-f'(q)(s-q)\le| f'(s)-f'(q)||s-q|\le f'(1)|s-q|.
	\end{align}

(a) We first consider the case $\lambda>\rho$.
          For any $\epsilon\in(0, (\lambda-\rho)\wedge(\rho/2)\wedge (\lambda/2\rho))$, let $t_\epsilon$ be the constant in Lemma \ref{lem:est-u}.
               Let $T>t_\epsilon$ be large enough so that $h(T)>x_0$.
By H\"older's inequality and \eqref{e:V},  we have for $t\ge t_0$ and $s\in (T, t-r(t))$,
\begin{align}
		&I_1(t,s):=\mE V(u(s+r(t),-c\,
		h(r(t))
		+\xi_{t-r(t)-s});|\xi_{t-r(t)-s}|>h(s+r(t))) \nonumber\\
		\le &f'(1) \mE(|u(s+r(t),-c
		h(r(t))
		+\xi_{t-r(t)-s})-q|;
		|\xi_{t-r(t)-s}|>h(s+r(t)))\nonumber\\
		\le & f'(1)\sqrt{\mE(|u(s+r(t),-c\,		h(r(t))
		+\xi_{t-r(t)-s})-q|^2)\mP(
		|\xi_{t-r(t)-s}|>h(s+r(t))
		)}\nonumber\\
		\le &f'(1)\sqrt{c_0c(\epsilon)s^{6\rho/\lambda}[e^{-2\rho(1-\epsilon)s}+t
		e^{-\lambda r(t)}
		e^{-(2\rho\wedge \lambda)(1-\epsilon)s}] t
		e^{-\lambda(s+r(t))}
		}\nonumber\\
		\le &f'(1)\sqrt{c_0c(\epsilon)t
		e^{-\lambda r(t)}}
		s^{3\rho/\lambda}\left[\sqrt{1+t
		e^{-\lambda r(t)}
		}e^{-\rho(1-\epsilon)s}e^{-\lambda s/2}+\sqrt{t
		e^{-\lambda r(t)}
		}e^{-\lambda(1-\epsilon/2)s}\right],
\end{align}
where in the third inequality we used Lemma \ref{lem:est-u} and \eqref{tail-prob}, and in the fourth inequality we used the inequality $\sqrt{a+b}\le \sqrt{a}+\sqrt{b}$,  $a,b\ge0$.
	Since $\sum_n n\Lambda(n)^{-\alpha}L(\Lambda(n))=\sum_n n e^{-\lambda r(n)}<\infty$ and $\Lambda$ is non-decreasing, we know that $te^{-\lambda r(t)}\to 0$ as $t\to\infty$.
	It follows that
	\begin{align}\label{small-int-1}
				&\int_T^{t-r(t)}
		\beta e^{\rho s}I_1(t,s)\,ds
		\le f'(1) \sqrt{c_0c(\epsilon)t
		e^{-\lambda r(t)}
		}\nonumber\\&\qquad \qquad \times \left[\sqrt{1+t
		e^{-\lambda r(t)}
		}\int_T^\infty \beta s^{3\rho/\lambda}e^{-(\lambda/2-\rho\epsilon)s}\,ds+\sqrt{t
		e^{-\lambda r(t)}
		}\int_T^\infty  \beta s^{3\rho/\lambda} e^{-(\lambda-\rho-\epsilon/2)s}\,ds\right]\nonumber\\
		&\to 0, \quad t\to\infty.
	\end{align}

	By \eqref{Rt/ht}, we have $\P(R_{t}\le 2h(t))\to \E(e^{-\vartheta^* 2^{-\alpha} W})\in(q,1)$.
	Thus there exists $a\in(q,1)$ such that when $T$ is large enough,
	 $$\P(R_{t}\le 2h(t))\le a, \quad t\ge T.$$
	If $|\xi_{t-r(t)-s}|\le h(s+r(t))$, then for any
	$t>t_0$ and $s\in  (T, t-r(t))$,
	\begin{align*}
	u(r(t)+s,-h(r(t))+\xi_{t-r(t)-s})&\le \P(R_{s+r(t)}\le h(r(t))+h(s+r(t)))\\
&\le \P(R_{s+r(t)}\le 2h(s+r(t)))\le a.
	\end{align*}
Thus by \eqref{est:V}, we have
	$$
	V(u(r(t)+s,-h(r(t))+\xi_{t-r(t)-s}))\le f''(a)(u(r(t)+s,-h(r(t))+\xi_{t-r(t)-s})-q)^2.
	$$
	By Lemma \ref{lem:est-u}, we have for any
	$t>t_0$ and $s\in  (T, t-r(t))$,
	\begin{align}
		I_2(t,s)&:=\mE(V(u(s+r(t),-c\,
		h(r(t))
		+\xi_{t-r(t)-s}));
		 |\xi_{t-r(t)-s}|\le
		h(s+r(t))
		)\nonumber\\
		&\le f''(a)\mE(u(s+r(t),-c\,
		h(r(t))
		+\xi_{t-r(t)-s})-q)^2\nonumber\\
		&\le f''(a)c(\epsilon)s^{6\rho/\lambda}[e^{-2\rho(1-\epsilon)s}+t
		 e^{-\lambda r(t)}
		e^{-(2\rho\wedge \lambda)(1-\epsilon)s}] \nonumber\\
		&\le f''(a)c(\epsilon)s^{6\rho/\lambda}[(1+t
		e^{-\lambda r(t)}
		)e^{-2\rho(1-\epsilon)s}+t
		e^{-\lambda r(t)}
		e^{- \lambda(1-\epsilon)s}].
	\end{align}
It follows that
	\begin{align}\label{small-int-2}
		&	
		\int_T^{t-r(t)}
		\beta e^{\rho s}I_2(t,s)\,ds\nonumber\\
		\le &  f''(a)c(\epsilon)\left[(1+t
		e^{-\lambda r(t)}
		)\int_T^\infty s^{6\rho/\lambda}e^{-(\rho-2\epsilon)s}\,ds+ t
		e^{-\lambda r(t)}
		 \int_T^\infty s^{6\rho/\lambda}e^{-(\lambda-\rho-\epsilon)s}\,ds\right]\nonumber\\
		\to & f''(a)c(\epsilon)\int_T^\infty s^{6\rho/\lambda}e^{-(\rho-2\epsilon)s}\,ds,\quad t\to\infty.
	\end{align}
By \eqref{small-int-1} and \eqref{small-int-2}, we have
\begin{align}\label{small-int}
	&\limsup_{t\to\infty}e^{\rho(t-r(t))}J_1(t, T)\le
	f''(a)c(\epsilon)\int_T^\infty s^{6\rho/\lambda}e^{-(\rho-2\epsilon)s}\,ds.
\end{align}

 (b) We now assume $\lambda=\rho$. Since $\Lambda(t)\overset{t}{>}e^{\gamma t}$, we have
 $$e^{-\lambda r(t)}=\Lambda(t)^{-\alpha}L(\Lambda(t))\overset{t}{\le }\Lambda(t)^{-\alpha/2}\le e^{-\alpha \gamma t/2},$$
 which implies that $r(t)\overset{t}{\ge }\frac{\alpha\gamma}{2\lambda}t$.
  Let $\epsilon\in(0, \frac{\alpha\gamma}{2\lambda}\wedge \frac{\rho}{2})$,
  and let $T>t_\epsilon$.

 By Lemma \ref{lem:lambda-rho}, $p_k=0$ for $k\ge 3$. Thus $f''(1)<\infty$. By \eqref{est:V},
 we have, for any $t\ge t_0$ and $s\in  (T, t-r(t))$,
 $$
 V(u(r(t)+s,-h(r(t))+\xi_{t-r(t)-s}))\le f''(1)(u(r(t)+s,-h(r(t))+\xi_{t-r(t)-s})-q)^2.
 $$
 Applying Lemma \ref{lem:est-u}, we have for any
 $t\ge t_0$ and $s\in (T, t-r(t))$,
 \begin{align}
 	&\mE(V(u(s+r(t),-c\,h(r)+\xi_{t-r(t)-s})))
 	\le f''(1)c(\epsilon)s^{6\rho/\lambda}[ e^{-2\rho(1-\epsilon)s}+t
	e^{-\lambda r(t)}
	e^{- \rho(1-\epsilon)s}].
 \end{align}
By \eqref{673}, it follows that
\begin{align}\label{small-int'}
	&e^{\rho(t-r(t))}J_1(t, T)\le
	f''(1)c(\epsilon)\left[\int_T^\infty s^{6\rho/\lambda}e^{-(\rho-2\epsilon)s}\,ds+ te^{-\lambda r(t)}
	 \int_T^{t-r(t)}s^{6\rho/\lambda}
	e^{\rho\epsilon s}\,ds\right]\nonumber\\
	\le &  f''(1)c(\epsilon)\left[\int_T^\infty s^{6\rho/\lambda}e^{-(\rho-2\epsilon)s}\,ds+ (\rho\epsilon)^{-1}t^{1+6\rho/\lambda}e^{-\alpha \gamma t/2}e^{\lambda\epsilon t}\right]\nonumber\\
	\to & f''(1)c(\epsilon)\int_T^\infty s^{6\rho/\lambda}e^{-(\rho-2\epsilon)s}\,ds,\quad t\to\infty.
\end{align}

	Combining \eqref{large-int},\eqref{mid-int},\eqref{small-int}  and \eqref{small-int'},
	and letting  $T\to\infty$, we obtain that
	\begin{align}\label{649}
		\lim_{t\to\infty}e^{\rho(t-r(t))}J(t)= \int_{-\infty}^\infty\beta e^{\rho s} V \left[\phi(C(\varphi)e^{\lambda s})\right]\,ds=A(\phi(C(\varphi))).
	\end{align}
Thus by \eqref{648}, \eqref{647}	 and \eqref{649}, we have that
$$\lim_{t\to\infty}e^{\rho (t-r(t))}\E^*\Big(\prod_{u\in\mathcal{L}_t}\varphi(\xi_t^u/\Lambda(t))\Big) = \frac{1}{1-q} A\left[\phi(C(\varphi))\right] .$$

Since $C(\varphi) $ is continuous in  $\varphi$, the above limit also valid for $\varphi(y)={\bf 1}_{(-\infty,1]}(y)$, that is,
$$\lim_{t\to\infty}e^{\rho (t-r(t))}\E^*\Big(R_t\le \Lambda(t)\Big) = \frac{1}{1-q} A(\phi(q_1\vartheta/\alpha)) .$$
The proof is now complete.
	
	\hfill $\Box$
	
	\noindent{\bf Proof of Corollary \ref{cor:lower deviation}:} 	For any $g\in C_c^+(\overline{\R}_0)$ and $\theta>0$, applying Theorem \ref{them:lower deviation}  with $\varphi(y)=e^{-\theta g(y)}{\bf 1}_{(-\infty, 1]}(y)$  we get
	\begin{align*}
		\lim_{t\to\infty}\E^*\left(e^{-\theta\frac{\mathbb{X}_t}{\Lambda(t)}(g)}|R_t\le \Lambda(t)\right)=	\lim_{t\to\infty}\frac{e^{\rho(t-r(t))}\E^*\left(e^{-\theta\frac{\mathbb{X}_t}{\Lambda(t)}(g)};R_t\le \Lambda(t)\right)}{e^{\rho(t-r(t))}\P^*(R_t\le \Lambda(t))}=\frac{A(\phi(C(\varphi)))}{A(\phi(\vartheta^*))}.
	\end{align*}
By the definition of $\Xi$, we have
	\begin{align}\label{650}
		E(e^{-\theta \Xi(g)})=\sum_{j=1}^\infty P(\mathcal{K}=j)\left[E(e^{-\theta
		\bar{\mathcal{N}}_\infty(g)})\right]^j=\frac{A\left((\phi(\vartheta^*)-q)E(e^{-\theta \bar{\mathcal{N}}_\infty(g)})+q\right)}{A(\phi(\vartheta^*))}.
		\end{align}
	Note that $\{\mathcal{N}_\infty(\R)>0\}=\{W>0\}$. Thus
	\begin{align*}
	E(e^{-\theta \bar{\mathcal{N}}_\infty(g)})&=E(e^{-\theta \mathcal{N}_\infty(g)}| \mathcal{N}_\infty(\R)>0, \mathcal{N}_\infty(1,\infty)=0)\\
	=&\frac{E(e^{-\theta \mathcal{N}_\infty(g)};\mathcal{N}_\infty(\R)>0, \mathcal{N}_\infty(1,\infty)=0)}{P(\mathcal{N}_\infty(\R)>0, \mathcal{N}_\infty(1,\infty)=0)}\\
	=&\frac{E(e^{-\theta \mathcal{N}_\infty(g)}; \mathcal{N}_\infty(1,\infty)=0)-P(\mathcal{N}_\infty(\R)=0)}{P( \mathcal{N}_\infty(1,\infty)=0)-P(\mathcal{N}_\infty(\R)=0)}\\
	=&\frac{\E(\exp\{-C(\varphi)W\})-\P(W=0)}{\P(\exp\{-\vartheta^* W\})-\P(W=0)}=\frac{\phi(C(\varphi))-q}{\phi(\vartheta^*)-q}.
	\end{align*}
By \eqref{650}, we have that
$$	E(e^{-\theta \Xi(g)})=\frac{A(\phi(C(\varphi)))}{A(\phi(\vartheta^*))}.$$
The proof is complete.
\hfill$\Box$

 \subsection{Lower deviation of $\Y_t$}

 For the proof of  Theorem \ref{them:large-deviation},
 we first established  the upper deviation for  $\mathbb{Y}_t$, and then used it to get the corresponding result for $\mathbb{X}_t$.
 However we can not get the lower  deviation  of $\X_t$ from that of $\Y_t$. In Section \ref{subsec3.1},  we proved Theorem \ref{main-theorem-inf}, the
lower large deviation result of  $\mathbb{X}_t$.
In this subsection,
we establish the lower deviation result of  $\mathbb{Y}_t$.
 Recall that $\phi$ is defined by \eqref{def-phi}.
\begin{proposition}
If $\Lambda:[0,\infty)\to(0,\infty)$ satisfies $\Lambda(t)\to\infty$ and
		$$t^2\Lambda(t)^{-\alpha}L(\Lambda(t))\to0,\quad \Lambda(t)/h(t)\to 0, $$  as $t\to\infty$,
then for any $\varphi\in\cH(\R)$,
		\begin{align*}
			\lim_{t\to\infty}e^{\rho(t-r(t)}\E^*( \cI(\varphi,\mathbb{Y}_t/\Lambda(t)))= A\left( \phi\left(C(\varphi)\right)\right) ,
		\end{align*}
	where $r=r(t)$ is defined by $\Lambda(t)=h(r(t)).$
In particular,
	$$\lim_{t\to\infty}e^{\rho (t-r(t))}
 \P^*\Big(M_t\le \Lambda(t)\Big) = A\left[\phi(\vartheta^*)\right].
$$
	\end{proposition}
	
	\noindent{\bf Proof:}
Since $\varphi\in\cH(\R)$, there exists $\delta>0$ such that $\varphi\equiv1$ on  $[-\delta, \delta].$
Note that
	\begin{align}
			(1-q)\E^*( \cI(\varphi,\mathbb{Y}_t/\Lambda(t)))=	\E( \cI(\varphi,\mathbb{Y}_t/\Lambda(t)))-q+q-	\E( \cI(\varphi,\mathbb{Y}_t/\Lambda(t));\mathcal{S}^c).
	\end{align}
On the event $\{Z_t=0\}$, we have $\cI(\varphi,\mathbb{Y}_t/\Lambda(t))=1$, thus
\begin{align}
	0\le q-	\E( \cI(\varphi,\mathbb{Y}_t/\Lambda(t));\mathcal{S}^c)\le
	q-\P(Z_t=0)\sim -A(0)e^{-\rho t}=
	o(e^{-\rho(t-r(t))}).
	\end{align}
By the definition of  $\mathbb Y_{r,t}$ given by  \eqref{314}, we have
	\begin{align}\label{e:YYprime}
	\E( \cI(\varphi,\mathbb{Y}_t/\Lambda(t)))\le \E( \cI(\varphi,
	 \mathbb{Y}_{r(t),t}/\Lambda(t)))
	&=\E\left( \left(\E \cI(\varphi,
	\mathbb{Y}'_{r(t)}/h(r(t)))\right)^{Z_{t-r(t)}}\right),
	\end{align}
	where $\mathbb{Y}'_{t}$ is defined in \eqref{def:Y-stprime}.
	If $Z_{r(t)}=0$, then $\cI(\varphi,\mathbb{Y}_{r(t)}/h(r(t)))=\cI(\varphi,\mathbb{Y}'_{r(t)}/h(r(t)))=1$; if $Z_{r(t)}>0$ then
	\begin{align*}
	\cI(\varphi,\mathbb{Y}_{r(t)}/h(r(t)))=&\cI(\varphi,\mathbb{Y}'_{r(t)}/h(r(t)))\cdot \varphi(X_{o,r(t)}/h(r(t)))^{Z_{r(t)}}\\
	\ge &\cI(\varphi,\mathbb{Y}'_{r(t)}/h(r(t)))\cdot{\bf 1}_{\{|X_{o,r(t)}|\le \delta h(r(t))\}}.
	\end{align*}
	By \eqref{tail-prob} and \eqref{def:h} we have
	$$
	 \P(|X_{o,r(t)}|>\delta h(r(t)))\overset{t}{\le}
c_0 r(t) (\delta h(r(t)))^{-\alpha}L(\delta h(r(t)))
\overset{t}{\le}
2c_0  \delta^{-\alpha}\cdot r(t)e^{-\lambda r(t)}\to0.$$
	Thus by \eqref{laplace-N3} we have
	\begin{align}\label{laplace-Y4}
	\lim_{t\to\infty}	
	\E(\cI(\varphi,
	 \mathbb{Y}'_{r(t)}/h(r(t))))&=\lim_{t\to\infty}\E(\cI(\varphi,\mathbb{Y}_{r(t)}/h(r(t))))=
	\phi(C(\varphi))<1.
	\end{align}
	Thus, by  \eqref{limit-F}, we have that
	\begin{align*}
		\lim_{t\to\infty}
		e^{\rho(t-r(t))}
		\left(\E\left( \left[\E \cI(\varphi,
		\mathbb{Y}'_{r(t)}/h(r(t)))\right]^{Z_{t-r(t)}} \right)-q\right)
		=A\left( \phi( C(\varphi))\right) .
	\end{align*}
Hence by \eqref{e:YYprime} we have
	\begin{align}
		\limsup_{t\to\infty}
		e^{\rho(t-r(t))}(
		\E( \cI(\varphi,\mathbb{Y}_t/\Lambda(t)))-q)\le A\left( \phi\left(C(\varphi)\right)\right) .
	\end{align}

	On the other hand,
\begin{align}\label{656}
		&\E(\cI(\varphi,\mathbb{Y}_t/\Lambda(t))|\FF_{t}^{\mathbb{T}})\nonumber\\
\ge &\P(\max_{u\in\mathcal{D}_t:b_u\le t-r(t)} |X_{u,t}|\le \delta h(r(t))|\FF_{t}^{\mathbb{T}})\E( I(\varphi,\mathbb{Y}_{r(t),t}/h(r(t)))|\FF_{t}^{\mathbb{T}}).
	\end{align}
Note that
 $\{u\in\mathcal{D}_t:b_u\le t-r(t)\}\subset \mathcal{D}_{t-r(t)}$. Thus
by \eqref{tail-prob2} we have that
	\begin{align*}
	\P(\max_{u\in\mathcal{D}_t:b_u\le t-r(t)} |X_{u,t}|\le \delta h(r(t))|\FF_{t}^{\mathbb{T}})\ge & \prod_{u\in\mathcal{D}_{t-r(t)}}\P(|X_{u,t}|\le\delta h(r(t)))|\FF_{t}^{\mathbb{T}})\\
 \overset{t}{\ge }&(1-2c_0 \delta^{-\alpha-1} t e^{-\lambda r(t)})^{|\mathcal{D}_{t-r(t)}|}.
	\end{align*}
By \eqref{656}, we have that
	\begin{align}
		&\E(\cI(\varphi,\mathbb{Y}_t/\Lambda(t)))
		\overset{t}{\ge}
		\E\left((1-2c_0\delta^{-\alpha}te^{-\lambda r(t)})^{|\mathcal{D}_{t-r(t)}|}\cdot  \cI(\varphi,\mathbb{Y}_{r(t),t}/\Lambda(t)) \right)\nonumber \\
		=&\E\left((1-2c_0\delta^{-\alpha}t
		e^{-\lambda r(t)})^{|\mathcal{D}_{t-r(t)}|}\cdot  [\E \cI(\varphi,\mathbb{Y}'_{r(t)}/\Lambda(t))]^{Z_{t-r(t)}} \right)\nonumber \\
			\ge& \E \left((1-2c_0\delta^{-\alpha}t
			e^{-\lambda r(t)})^{c\cdot(t-r(t)) Z_{t-r(t)}} [\E \cI(\varphi,\mathbb{Y}'_{r(t)}/\Lambda(t))]^{Z_{t-r(t)}} ;|\mathcal{D}_{t-r(t)}| \le c\cdot(t-r(t)) Z_{t-r(t)} \right) \nonumber\\
		\ge &
		\E\left(C(t)^{Z_{t-r(t)}} \right) -\P(|\mathcal{D}_{t-r(t)}|> c\cdot (t-r(t)) Z_{t-r(t)}),
	\end{align}
	where
	$C(t)=(1-2c_0\delta^{-\alpha}te^{-\lambda r(t)})^{c\cdot(t-r(t)) }\E \cI(\varphi,\mathbb{Y}'_{r(t)}/\Lambda(t))$  and $c>0$ is a constant.
	Since $t^2 e^{-\lambda r(t)}=t^2\Lambda(t)^{-\alpha}L(\Lambda(t))\to 0,$
by \eqref{laplace-Y4}, we have
	$$\lim_{t\to\infty}C(t)=\lim_{t\to\infty}\E \cI(\varphi,\mathbb{Y}'_{r}/\Lambda(t))=\phi\left(C(\varphi)\right) <1.$$
	By \eqref{limit-F}, we have that
	$$
	\lim_{t\to\infty}e^{\rho(t-r(t))}(\E\left( C(t)^{Z_{t-r(t)}} \right)-q)=A\left( \phi\left(C(\varphi)\right)\right).
	$$
	 Note that $|\mathcal{D}_{t-r(t)}|\le \sum_{u\in\mathcal{L}_{t-r(t)}}|I_u|$.
	Thus by the many-to-one formula,
	\begin{align*}
		&\P(|\mathcal{D}_{t-r(t)}|> c\cdot (t-r(t)) Z_{t-r(t)})\le \P\left(\sum_{u\in\mathcal{L}_{t-r(t)}}|I_u|>c\cdot (t-r(t))Z_{t-r(t)}\right)
		\\&\le
             \P\left(\sum_{u\in\mathcal{L}_{t-r(t)}}{\bf 1}_{\{|I_u|>c\cdot (t-r(t))\}}\ge 1\right)\le \E\left(\sum_{u\in\mathcal{L}_{t-r(t)}}{\bf 1}_{\{|I_u|>c\cdot (t-r(t))\}}\right)\\
		&\le
		 e^{\lambda(t-r(t))}P(1+n_{t-r(t)}>c \cdot (t-r(t)))\\
		&\le
		e^{\lambda(t-r(t))} e^{-c\cdot (t-r(t))}E(e^{(1+n_{t-r(t)})})
		=e\cdot e^{\lambda(t-r(t))}e^{-c\cdot(t-r(t))}e^{(e-1)\beta(t-r(t))}\\
		&
		=e\cdot e^{-(c-\lambda-(e-1)\beta)(t-r(t))}.
	\end{align*}
	Now we choose $c$ such that $c>\lambda+(e-1)\beta-\rho.$ Then
	$$
		 \lim_{t\to\infty}e^{\rho(t-r(t))}\P(|\mathcal{D}_{t-r(t)}|>c(t-r(t)) Z_{t-r(t)})=0.
	$$
Hence we have
	\begin{align}
		\liminf_{t\to\infty}
		e^{\rho(t-r(t))}
		(\E(\cI(\varphi,\mathbb{Y}_t/\Lambda(t)))-q)\ge A\left( \phi\left(C(\varphi) \right)\right).
	\end{align}

The proof is complete now.

\hfill$\Box$
	
	\section{Almost  sure convergence results}
	In this section we will give the proofs of Theorems \ref{main-theorem-inf} and \ref{them:upper-limit}. We first prove the corresponding results for $\mathbb{Y}$.
	Note that for any $a>0$,
	$$\{R_t\le a\}=\cI({\bf 1}_{(-\infty, a]},\X_t), \quad \{M_t\le a\}=\cI({\bf 1}_{(-\infty, a]},\Y_t).$$
	By the continuity of $C(\varphi)$, the results established in the previous sections for $\cH(\R)$ are also valid for ${\bf 1}_{(-\infty, a]}$.
	
	\subsection{Almost sure convergence results for $M_t$}
		\begin{lemma}\label{prop:upper-limit}
		Suppose  that $G:[0,\infty)\to (0,\infty)$
		satisfies
		$\lim_{t\to\infty}\frac{G(t)}{h(t)}=\infty$.
		\begin{itemize}
			\item[(1)] If $\sum_{n}e^{\lambda n}G(n)^{-\alpha}L(G(n))<\infty$, then
			\begin{align*}
				\limsup\limits_{t \to \infty}\dfrac{ M_n}{G(n)}\le 0,  \qquad \P \textup{-a.s.}
			\end{align*}
			\item[(2)] If  $\sum_{n}e^{\lambda n}G(n)^{-\alpha}L(G(n))=\infty$, then
			\begin{align*}
				\limsup\limits_{t \to \infty}\dfrac{ M_n}{G(n)}= \infty,  \qquad \P^* \textup{-a.s.}
			\end{align*}
			
		\end{itemize}
		\end{lemma}
	
	\noindent{\bf Proof:}
	(1) Assume that $\sum_{n}e^{\lambda n}G(n)^{-\alpha}L(G(n))<\infty$.
		By Proposition \ref{prop:large-M}, we have  for any $c>0$,
	$$\lim_{t\to\infty}e^{-\lambda t}G(t)^{\alpha}L(G(t))^{-1}
	 \P(M_t>cG(t))
	=c^{-\alpha}\vartheta^*,
$$
where $\vartheta^*$ is
given in \eqref{d:varthetastar}.
Combining this with $\sum_{n}e^{\lambda n}G(n)^{-\alpha}L(G(n))<\infty$, we get that for any $c>0$,
	$$\sum_{n=1}^\infty
	 \P(M_n>cG(n))<\infty.
	$$
	By the Borel-Cantelli lemma, we get that
	$$\limsup\limits_{n \to \infty}\dfrac{ M_n}{G(n)}\le c,  \qquad \P \textup{-a.s.}$$
	Letting $c\to 0,$ we get the desired result.
	
	(2)  Assume that $\sum_{n}e^{\lambda n}G(n)^{-\alpha}L(G(n))=\infty$. For $0<s<t$, define $$M_{s,t}^{(2)}:=\max_{u\in\cD_t: t-s<b_u\le t} X_{u,t}.$$ For any $c>0$, we have that
	\begin{align}\label{3.12}
		\P(M_n>cG(n)|\mathcal{F}_{n-1})\ge \P(M_{1,n}^{(2)}>cG(n)|\mathcal{F}_{n-1}).
	\end{align}
	By \eqref{limit-M2} with $s=1$, $t=n$, $\Lambda(t)=cG(t)$ and $\varphi={\bf 1}_{(-\infty,1]}$, we get that, as $n\to\infty$,
	\begin{align}\label{3.13}
		&e^{-\lambda n}G(n)^{\alpha}L^{-1}(G(n))\P(M_{1,n}^{(2)}>cG(n)|\mathcal{F}_{n-1})\nonumber\\
		\to&  c^{-\alpha}\dfrac{q_1}{\alpha} e^{-\lambda} \int_{0}^{1}  \left(e^{\lambda r}-e^{-\beta r}\right)\P(Z_{1-r}>0)\d r\cdot W\,\quad  \mbox{a.s.}.
	\end{align}
	Combining \eqref{3.12} and \eqref{3.13},  we obtain that, for any $c>0$, on the event $\{W>0\}$,
	\begin{align}\label{claim1}
		\sum_{n=1}^\infty \P(M_n>cG(n)|\mathcal{F}_{n-1})=\infty,\quad  \P\mbox{-a.s.}
	\end{align}
	Applying the conditional Borel-Cantelli lemma, we get that for any $c>0$,
	$$\P^*(M_n>cG(n),i.o.)=1.$$
Since $c$ is arbitrary, desired result follows immediately.
	
	\hfill$\Box$

	\begin{lemma}\label{lem:lower-inf} It holds that
		\begin{align*}
			\liminf\limits_{t \to \infty}\dfrac{M_n}{H(e^{-\lambda n}\log n)} \le \left(\vartheta^* W\right)^{\frac{1}{\alpha}}, \qquad \P^* \textup{-a.s.}
		\end{align*}
	\end{lemma}
\noindent{\bf Proof:}
For any $\epsilon\in(0,1)$, choose $\epsilon'>0$ small enough
such that
$(1+\epsilon^{\prime})(1+\epsilon/2)^{-\alpha}<1$.
Let $n_k:=\lfloor k^{1+\epsilon^{\prime}}\rfloor$, $k\ge 1$, where $\lfloor x\rfloor$ is the integer part of $x$. It is easy to see $n_k<n_{k+1}$ for any $k\ge 1$. Let
$$U_{k}:=H\left(\dfrac{\log n_k}{\vartheta^* Z_{n_{k-1}}e^{\lambda(n_k-n_{k-1})}}\right).$$
Since $H(y)=y^{-1/\alpha}\bar{L}(y^{-1})$,  we have
\begin{align}\label{3.112}
	\lim_{k\to\infty}\frac{U_{k}}{H(e^{-\lambda n_k}\log n_k)}=\lim_{k\to\infty}(\vartheta^* W_{n_{k-1}})^{1/\alpha}=(\vartheta^* W)^{1/\alpha},\qquad \P^*\mbox{-a.s.}
\end{align}	
			Note that $U_k\in\mathcal{F}_{n_{k-1}}$. Thus we have
	\begin{align}\label{Mnk}
		\P\left(M_{n_k}\leq (1+\epsilon)U_{k}|\FF_{n_{k-1} }\right)
		&\geq \1_{\{M_{n_{k-1}}\leq (1+\epsilon/2)U_{k}\}} J_k
	\end{align}
	where $$J_k:= \P\left(\forall v \in \LL_{n _{k-1}},X_{v,n_{k}}-X_{v,n_{k-1}} \leq
\epsilon U_{k}/2,
M^{v}_{n_k}\leq (1+\epsilon/2)U_{k} |\FF_{n_{k-1} }\right)$$ and $M^v_{t}=\max_{u\in \cD_t: v<u}X_{u,t}$.
	
	We first show that, on the survival event $\S$,
	$\{M_{n_{k-1}}> (1+\epsilon/2)U_{k}\}$ can occur only finitely many times. By \eqref{est-Zn}, we have for $p>\lambda/\rho$,
	\begin{align}\label{453-1}
		\sum_{n=1}^\infty
		\P(0<e^{-\lambda n} Z_n<n^{-p})
		\le C \sum_{n=1}^\infty n^{-p\rho/\lambda}<\infty.
	\end{align}
	Note that
	\begin{align*}
		&\P\left(M_{n_{k-1}}> (1+\epsilon/2)U_{k}, \S\right) \\
		&\leq
		\P\left(e^{-\lambda n_{k-1}}Z_{n_{k-1}}<n_{k-1}^{-p},
		\S\right)+\P\left(M_{n_{k-1}}> (1+\epsilon/2)U_k,
		e^{-\lambda n_{k-1}}Z_{n_{k-1}}\geq n_{k-1}^{-p}\right)\\
		& \leq \P\left(
		e^{-\lambda n_{k-1}}Z_{n_{k-1}}<n_{k-1}^{-p},\S\right)
		+\P\left(M_{n_{k-1}}>  (1+\epsilon/2)H\left(\dfrac{\log n_k}{\vartheta^* n_{k-1}^{-p}e^{\lambda n_k}}\right)\right).
	\end{align*}
	By \eqref{453-1}, we have
	\begin{align}\label{453}
		\sum\limits_{k=1}^{\infty}\P\left(
		e^{-\lambda n_{k-1}}Z_{n_{k-1}}<n_{k-1}^{-p},\S\right)
		\le \sum\limits_{k=1}^{\infty}\P\left(0<
		e^{-\lambda n_{k-1}}Z_{n_{k-1}}<n_{k-1}^{-p}\right)
		< \infty.
	\end{align}
	Observe that $n_k-n_{k-1}\sim (1+\epsilon')k^{\epsilon'}$. Thus
	$$\lim_{k\to\infty}\frac{H\left(\dfrac{\log n_k}{\vartheta^* n_{k-1}^{-p}e^{\lambda n_k}}\right)}{h( n_{k-1})}=\lim_{k\to\infty}\left(\frac{\dfrac{\log n_k}{\vartheta^* n_{k-1}^{-p}e^{\lambda n_k}}}{e^{-\lambda n_{k-1}}} \right)^{-1/\alpha} =\infty.$$
Hence by Proposition \ref{prop:large-M},
	we have that as $k \to \infty$,
	\begin{align}\label{M_{n_{k-1}}}
		&\P\left(M_{n_{k-1}}> (1+\epsilon/2)H\left(\dfrac{\log n_k}{\vartheta^* n_{k-1}^{-p}e^{\lambda n_k}}\right)\right)  \notag \\
		&\sim \vartheta^* e^{\lambda{n_{k-1}}}(1+\epsilon/2)^{-\alpha}H\left(\dfrac{\log n_k}{\vartheta^* n_{k-1}^{-p}e^{\lambda n_k}}\right)^{-\alpha}L\left(H\left(\dfrac{\log n_k}{\vartheta^* n_{k-1}^{-p}e^{\lambda n_k}}\right)\right) \notag \\
		& =(1+\epsilon/2)^{-\alpha}\dfrac{n_{k-1}^{p}\log{n_k}}{e^{\lambda(n_k-n_{k-1})}},
	\end{align}
where in the last equality we used \eqref{H-inverse}.
	Since $n_k=\lfloor k^{1+\epsilon^{\prime}}\rfloor$, we have $\sum_k\dfrac{n_{k-1}^{p}\log{n_k}}{e^{\lambda(n_k-n_{k-1})}}<\infty.$
	Combining \eqref{453} and \eqref{M_{n_{k-1}}}, we get that
	\begin{align*}
		\sum\limits_{k=1}^{\infty}\P\left(\{M_{n_{k-1}}> (1+\epsilon/2)U_{k}\}\cap\S\right)< \infty.
	\end{align*}
	Using the Borel-Cantelli lemma, we get that
	\begin{align}\label{Ak}
		\sum\limits_{k=1}^{\infty} \1_{\{M_{n_{k-1}}> (1+\epsilon/2)U_{k}\}} <\infty,  \quad \P\text{-a.s. on } \mathcal{S}.
	\end{align}

	Now we consider $J_k$. Recall that $M'_t$ is defined in \eqref{M_t'}.
	By the Markov property and the branching property, we have that
	\begin{align} \label{454}
		J_k
		=& \left[\P\left(X_{o,n_{k}-n_{k-1}} \leq\epsilon x/2,M'_{n_k-n_{k-1}}\leq (1+\epsilon/2)x \right)\right]^{Z_{n_{k-1}}}|_{x=U_{k}}\notag  \\
		=&e^{Z_{n_{k-1}}\log \P\left(X_{o,n_{k}-n_{k-1}} \leq\epsilon x/2,M'_{n_k-n_{k-1}}\leq (1+\epsilon/2)x \right)}|_{x=U_{k}} \notag \\
		\geq &e^{Z_{n_{k-1}}\log\left(1-\P\left(X_{o,n_{k}-n_{k-1}} >\epsilon x/2 \right)-\P\left(M'_{n_k-n_{k-1}}> (1+\epsilon/2)x\right)\right)}|_{x=U_{k}}.
	\end{align}
Note that,  on the survival event $\mathcal{S}$, $U_k\to\infty$.
Then by \eqref{tail-prob2},
we have for $k$ large enough,
	\begin{align}\label{456}
		\P\left(X_{o,n_{k}-n_{k-1}} >\epsilon x/2 \right)|_{x=U_{k}}\le& 2c_0 	(\epsilon/2)^{-\alpha-1} (n_{k}-n_{k-1})U_{k}^{-\alpha}L(U_{k})\nonumber\\
		=&2c_0 	(\epsilon/2)^{-\alpha-1} (n_{k}-n_{k-1})\dfrac{\log n_k}{\vartheta^* Z_{n_{k-1}}e^{\lambda(n_k-n_{k-1})}} \nonumber\\
		\to &0, \quad
		 \P\text{-a.s. on } \mathcal{S}
	\end{align}
as $k\to\infty$.
	By  Proposition \ref{prop:large-M} we have
	\begin{align}\label{457}
		\P\left(M^{\prime}_{n_k-n_{k-1}}> (1+\epsilon/2)x\right)|_{x=U_{k}}&\le \P\left(M_{n_k-n_{k-1}}> (1+\epsilon/2)x\right)|_{x=U_{k}}  \nonumber\\
		&\sim   \vartheta^* e^{\lambda(n_k-n_{k-1})} (1+\epsilon/2)^{-\alpha}U_{k}^{-\alpha}L(U_{k} ) \notag\\
		& = (1+\epsilon/2)^{-\alpha}\frac{\log n_k}
		{Z_{n_{k-1}}}\to 0, \quad \P\text{-a.s. on } \mathcal{S}
	\end{align}
	as $k\to\infty$.
	Since $n_k= \lfloor k^{1+\epsilon'}\rfloor$, using \eqref{456} and \eqref{457}, we obtain that
	\begin{align}\label{458}
		&\limsup_{k\to\infty}-\frac{1}{\log k}Z_{n_{k-1}}\log\left(1-\P\left(X_{o,n_{k}-n_{k-1}} >\epsilon x/2 \right)-\P\left(M^{\prime}_{n_k-n_{k-1}}> (1+\epsilon/2)x\right)\right)_{x=U_k}\nonumber\\
		=& \limsup_{k\to\infty}\frac{1}{\log k}Z_{n_{k-1}}\left(\P\left(X_{o,n_{k}-n_{k-1}} >\epsilon x/2 \right)+\P\left(M^{\prime}_{n_k-n_{k-1}}> (1+\epsilon/2)x\right)\right)_{x=U_k}\nonumber\\
		\le& (1+\epsilon')(1+\epsilon/2)^{-\alpha}<1, \quad
		\P\text{-a.s. on } \mathcal{S}.
	\end{align}
	Combining \eqref{454} and \eqref{458}, we get that,
	\begin{align}\label{459}
		\sum_{k=1}^\infty J_k=\infty, \quad
	 \P\text{-a.s. on } \mathcal{S}.
	\end{align}

	Combining \eqref{Mnk},\eqref{Ak} and \eqref{459},  we get that
	\begin{align*}
		\sum_{k}	\P\left(M_{n_k}\leq (1+\epsilon)U_{k}|\FF_{n_{k-1} }\right)=\infty, \quad \P\text{-a.s. on } \mathcal{S}.
	\end{align*}
	Applying the Borel-Cantelli lemma and then letting $\epsilon\to0$, we obtain that
	\begin{align*}
		\liminf_{k\to\infty}\frac{M_{n_k}}{U_k}\le 1\quad \P^*\text{-a.s.}
	\end{align*}
Hence the desired result follows from \eqref{3.112} immediately.

\hfill$\Box$

	We now modify the definition of $\mathcal{D}_t$ and $\mathcal{L}_t$ slightly.
For any $\delta\ge 0,$ and $t>0$, define
$$\mathcal{L}_{t,\delta}:=\{ u\in\mathcal{L}_{t}: Z^u_{t+\delta}>0\},$$
$$\mathcal{D}_{t,\delta}:=\{ u\in\mathbb{T}: b_u<t,  Z^u_{t+\delta}>0\}=\cup_{u\in \mathcal{L}_{t,\delta}} I_u,$$
$$R_{t,\delta}:=\max_{u\in \mathcal{L}_{t,\delta}}\xi_{t}^u,\quad  M_{t,\delta}:=\max_{u\in \mathcal{D}_{t,\delta}}X_{u,t}, \mbox{ and }
M'_{t,\delta}:=\max_{u\in \mathcal{D}_{t,\delta}\setminus\{o\}}X_{u,t}.
$$

Using argument similar to that in the proof of Lemma \ref{lem:Ms} with
$\varphi={\bf 1}_{(-\infty,1]}$,
we get the following result. We omit the proof.

\begin{lemma}\label{lem:Ms'}
	For any $\delta\ge 0$, it holds that
	$$\lim_{x\to\infty}x^{\alpha}L(x)^{-1}\P(M'_{s, \delta}>x)=\frac{q_1}{\alpha}\int_{0}^{s}  \left(e^{\lambda r}-e^{-\beta r}\right)\P(Z_{s+\delta-r}>0)dr
	=:c(s,\delta).
	$$
\end{lemma}

\begin{lemma} \label{lemma:low-Rn} For any $\delta\ge0$, it holds that
	\begin{align*}
		\liminf\limits_{n \to \infty}\dfrac{M_{n\delta,\delta }}{H(e^{-\lambda n\delta}\log n)} \geq  \left(\vartheta_\delta W\right)^{\frac{1}{\alpha}} \qquad \P^* \textup{-a.s.},
	\end{align*}
	where $\vartheta_\delta=\frac{q_1}{\alpha}\int_0^\infty e^{-\lambda r}\P(Z_{r+\delta}>0)\,dr.$
\end{lemma}
\noindent{\bf Proof:}
For  $1\le m<n$, let $U_{n,m}= H\left(\dfrac{\log n}{Z_{(n-m)\delta}c(m\delta,\delta) }\right)$ when $Z_{(n-m)\delta}>0$ and $U_{n,m}=0$ when $Z_{(n-m)\delta}=0$.  Note that
$$M_{n\delta,\delta}\ge \max_{u\in\mathcal{L}_{(n-m)\delta}}\max_{v\in \mathcal{D}_{n\delta,\delta}^u }X_{v,n\delta},$$
where $\mathcal{D}^u_{t,\delta}=\{v\in\mathcal{D}_{t,\delta:}: v>u\}.$
By the Markov property and the branching property, we have that,
conditioned on $\mathcal{F}_{(n-m)\delta}$,  $\{\max_{v\in \mathcal{D}_{n\delta,\delta}^u }X_{v,n\delta}, u\in\mathcal{L}_{(n-m)\delta}\}$ are i.i.d. with the same law as
$M'_{m\delta,\delta}$. Let $G_{m,\delta}(x) := \P(M_{m\delta,\delta}^{'}>x)$. Note that $U_{n,m}\in \FF_{(n-m)\delta}$. Thus for any $\epsilon\in(0,1)$, we have
\begin{align}\label{417}
	&\P\left(M_{n\delta,\delta}<(1-\epsilon)U_{n,m}|\FF_{(n-m)\delta} \right)
	 \le  \left(1-G_{m,\delta}((1-\epsilon)U_{n,m})\right)^{Z_{(n-m)\delta}}\nonumber\\
	\leq & e^{-Z_{(n-m)\delta} G_{m,\delta}((1-\epsilon)U_{n,m})}.
\end{align}
Note that on the survival event $\mathcal{S}$, $\lim_{n\to\infty}U_{n,m}=\infty$. Thus by Lemma \ref{lem:Ms'}, we have that, on the survival event $\mathcal{S}$, as $n\to\infty$,
$$Z_{(n-m)\delta} G_{m,\delta}((1-\epsilon)U_{n,m})\sim Z_{(n-m)\delta} c(m\delta,\delta) (1-\epsilon)^{-\alpha}U_{n,m}^{-\alpha}L(U_{n,m})=(1-\epsilon)^{-\alpha}\log n,\quad  \P\mbox{-a.s.}$$
where in the last equality we used \eqref{H-inverse}.
Since $(1-\epsilon)^{-\alpha}>1$,
we have that, on the survival event $\mathcal{S}$,   for any $m\ge1$, $$\sum_{n=m+1}^\infty\P\left(M_{n\delta,\delta}<(1-\epsilon)U_{n,m}|\FF_{(n-m)\delta}\right)\le \sum_{n=m+1}^\infty e^{-Z_{(n-m)\delta} G_{m,\delta}((1-\epsilon)U_{n,m})}<\infty, \quad  \P\mbox{-a.s.}.$$
By the second Borel-Cantelli lemma, we have that,  on the survival event $\mathcal{S}$,
\begin{align}\label{3.43}
	\sum\limits_{n=m+1}^{\infty}\1_{\{M_{n\delta,\delta}<(1-\epsilon)U_{n,m}\}}< \infty, \quad \P\mbox{-a.s.}
\end{align}
Thus,  by the definition of $H$, we obtain that under $\P^*$, for any $\epsilon\in(0,1)$ and  integer $m\ge1$,
\begin{align}\label{3.222}
	\liminf\limits_{n \to \infty}\dfrac{M_{n\delta,\delta} }{H(e^{-\lambda n\delta}\log n)}&\ge \liminf\limits_{n \to \infty}\dfrac{(1-\epsilon)U_{n,m} }{H(e^{-\lambda n\delta}\log n)}=(1-\epsilon)\lim_{n\to\infty}\left( \frac{\dfrac{\log n}{Z_{(n-m)\delta}c(m\delta,\delta)}}{e^{-\lambda n\delta}\log n}\right)^{-1/\alpha}\nonumber\\
	&=(1-\epsilon)\left( We^{-\lambda m\delta}c(m\delta,\delta)\right) ^{1/\alpha}, \quad \P^*\mbox{-a.s.}
\end{align}
Note that $\lim_{m\to\infty}e^{-\lambda m\delta}c(m\delta,\delta)=\vartheta_\delta.$		
Letting $m\to\infty$ and $\epsilon\to0$ in \eqref{3.222}, we get the desired result.
The proof is now complete.

\hfill$\Box$

	\subsection{Proofs of   Theorems \ref{main-theorem-inf} and  \ref{them:upper-limit}}
	
		\begin{lemma}\label{cor:M-R}
		If $\{a_n: n\ge 1\}$ is a sequence of positive numbers satisfying
		$\sum_{n}e^{\lambda n} a_n^{-p\alpha}<\infty$ for some $p\in (0, 2)$, then for any $\delta\ge 0$,
		\begin{align}\label{R-M}
			\limsup\limits_{n\to \infty}\dfrac{R_{n,\delta}-M_{n,\delta}}{a_n} \leq 0,\quad \P^*\mbox{-a.s.}
		\end{align}
		and
		\begin{align}\label{M-R}
			\limsup\limits_{n\to \infty}\dfrac{M_{n,\delta}-R_{n,\delta}^+}{a_n} \leq 0, \quad \P^*\mbox{-a.s.}
		\end{align}
	\end{lemma}
	\noindent{\bf Proof:}
	 Let $a(t)$ be a positive function satisfying $a(t)\overset{t}{\ge }e^{\epsilon t}$ for some $\epsilon\in (0, 1)$.
	 Let $c>2\alpha+e^2\beta/\epsilon^2$ be a constant. We choose $\theta\in (0, 1)$ so that  $\theta c< 1$.  We define the events $A_t(\theta)$ and $B_t$ as in the proof of Lemma \ref{lem:one-big'}.
	 We claim that for $n$ sufficiently large, on the event $A_t(\theta) \cap B_t\cap\{Z_t\ge 1\}$,
	 $$\{R_{t,\delta}-M_{t,\delta}>a(t)\}\cup\{M_{t,\delta}-R_{t,\delta}^+>a(t)\}= \emptyset.$$
	 In fact, by \eqref{35},  we have that on the event $A_t(\theta) \cap B_t\cap \{Z_t\ge1\}$,
	 \begin{align}\label{36}
	 	R_{t,\delta}=\max_{v\in\LL_{t,\delta}}\xi_t^v\le \max_{v\in\LL_{t,\delta}}X_{v^{\prime},t}+a(t)\leq   M_{t,\delta} +a(t),
	 \end{align}	
	 and
	 \begin{align}\label{37}
	 	M_{t,\delta} &= \max\limits_{v \in \LL_{t,\delta}} \max\limits_{u \in I_v} X_{u,t} =\left(\max\limits_{v \in \LL_{t,\delta}} X_{v^{\prime},t}\right) \vee \left(\max\limits_{v \in \LL_{t,\delta},u\in I_v \setminus\{v^{\prime}\}}X_{u,t}\right) \notag\\
	 	& \leq \left(\max\limits_{v \in \LL_{t,\delta}} \xi_t^v+a(t) \right)\vee \left(\theta a(t)/\log a(t)\right)\le R_{t,\delta}^++a(t),
	 \end{align}
	 for $t$ large enough so that $\log a(t)>\theta$.
	 Hence, by \eqref{36} and \eqref{37}, the claim is true.
	
	 Combining \eqref{31-1},  \eqref{Bto0}, we get that	for any $p\in(0,2)$,	
	 \begin{align}
	 	&\limsup\limits_{t \to \infty} e^{-\lambda t}a(t)^{p\alpha} \P\left(\left(\{R_{t,\delta}-M_{t,\delta}>a(t)\}\cup\{M_{t,\delta}-R_{t,\delta}^+>a(t)\}\right)\cap\{Z_t\ge1\}\right) \notag\\
	 	\leq& \limsup\limits_{t \to \infty}e^{-\lambda t}a(t)^{p\alpha}\{\P(A_t(\theta)^c) + \P(B_t^c) \} =0. \notag
	 \end{align}
	Thus we have
	\begin{align}
		\lim\limits_{t \to \infty} e^{-\lambda t}a(t)^{p\alpha} \P\left(\left\{R_{t,\delta}-M_{t,\delta}>a(t)\}\cup\{M_{t,\delta}-R_{t,\delta}^+>a(t)\}\right)\cap\{Z_t\ge1\}\right)=0. \notag
	\end{align}
Since $\sum_ne^{\lambda n} a_n^{-  p\alpha}<\infty, $ we have $a_n\overset{n}{\ge }e^{\frac{\lambda n}{p\alpha}}$. It follows that for any $\epsilon\in(0, 1)$ and $p\in(0,2)$,
	$$\P^*\left(R_{n,\delta}-M_{n,\delta}>\epsilon a_n\right) \le \frac{1}{\P(\mathcal{S})}\P\left( R_{n,\delta}-M_{n,\delta}>\epsilon a_n,Z_n\ge1\right) =o(e^{\lambda n} a_n^{-  p\alpha}).$$
	So we have that
	$$\sum_{n=1}^\infty\P^*\left(R_{n,\delta}-M_{n,\delta}>\epsilon a_n\right)<\infty. $$
	Applying the Borel-Cantelli lemma and then letting $\epsilon\to0$, we get  \eqref{R-M}.
	Similarly, we can get \eqref{M-R} as well.
	
	\hfill$\Box$

\bigskip

	\noindent {\bf Proof of Theorem \ref{main-theorem-inf}:} Note that $$H(e^{-\lambda n\delta}\log n)=e^{\lambda n\delta/\alpha}(\log n)^{-1/\alpha}\bar{L}(e^{\lambda n\delta}/\log n).$$
		Applying Lemma \ref{cor:M-R}  with $\delta=0$ and $\delta\ge 0$, we have
		\begin{align}\label{618}
\limsup_{n\to\infty}\frac{R_{n}-M_{n}}{H(e^{-\lambda n}\log n)}\le 0,\qquad		\limsup_{n\to\infty}\frac{M_{n\delta,\delta}-R^+_{n\delta,\delta}}{H(e^{-\lambda n\delta}\log n)}\le 0, \qquad \P^* \textup{-a.s.}
		\end{align}
By Lemma \ref{lem:lower-inf} we have
$$\liminf_{t\to\infty}\frac{R(t)}{H(e^{-\lambda t}\log t)}\le \liminf_{t\to\infty}\frac{M_n}{H(e^{-\lambda n}\log n)}+ \limsup_{n\to\infty}\frac{R_{n}-M_{n}}{H(e^{-\lambda n}\log n)}\le (\vartheta^* W)^{1/\alpha}, \quad \P^* \textup{-a.s.}$$

We now consider the lower bound of 		$\liminf_{t\to\infty}\frac{R(t)}{H(e^{-\lambda t}\log t)}$.
By Lemma \ref{lemma:low-Rn}  we have that  for any $\delta>0,$
\begin{align}\label{678}
	\liminf_{n\to\infty}\frac{R^+_{n\delta,\delta}}{H(e^{-\lambda n\delta}\log n)}\ge 	\liminf_{n\to\infty}\frac{M_{n\delta,\delta}}{H(e^{-\lambda n\delta}\log n)}-	\limsup_{n\to\infty}\frac{M_{n\delta,\delta}-R^+_{n\delta,\delta}}{H(e^{-\lambda n\delta}\log n)}\ge (\vartheta_\delta W)^{1/\alpha},\quad \P^* \textup{-a.s.}
	\end{align}
	For any $n\ge1$,  let $ Q_{n,\delta} =\inf_{n\delta\leq t \leq (n+1)\delta}R_t$.
	We claim  that  for any  $\epsilon\in(0,1)$ and $\delta>0$,
		\begin{align}\label{sum-Q_n}
			\sum\limits_{n=1}^{\infty} \P\left(Q_{n,\delta}<(1-\epsilon)H(e^{-\lambda n\delta}\log n)\left(\vartheta_\delta W_{n\delta}\right)^{\frac{1}{\alpha}}-H(n^{-2}) |\FF_{n\delta }\right)<\infty,\qquad \P^* \textup{-a.s.}
		\end{align}
Let $B_n:= \{ R_{n\delta,\delta}<(1-\epsilon)H(e^{-\lambda \delta n}\log n)\left(\vartheta_\delta W_{n\delta}\right)^{\frac{1}{\alpha}}\}$.	Note that
		\begin{align}\label{419}
			&\P\left(Q_{n,\delta}<(1-\epsilon)H(e^{-\lambda n\delta}\log n)\left(\vartheta_\delta W_{n\delta}\right)^{\frac{1}{\alpha}}-H(n^{-2}) |\FF_{n \delta}\right) \notag \\
			 \leq& \P(B_n|\FF_{n \delta})+ \P\left(B_n^c, Q_{n,\delta}<(1-\epsilon)H(e^{-\lambda n\delta}\log n)\left(\vartheta_\delta W_{n\delta}\right)^{\frac{1}{\alpha}}-H(n^{-2}) |\FF_{n\delta } \right).
		\end{align}
		By \eqref{678} and the Borel-Cantelli lemma, we know that on $\mathcal{S}$,
		\begin{align}\label{620}
			\sum\limits_{n=1}^{\infty}\P(B_n|\FF_{n \delta})<\infty, \qquad \P \textup{-a.s.}
		\end{align}
Now we consider the second term of the right hand side of  \eqref{419}.
Take  $Z_0=\nu=\sum_{k=1}^K\delta_{a_k}$ with
$a_1\ge a_2\ge\cdots \ge a_K$.
We use $Z^{(j)}_t$ to denote  the number of descendants at time $t$ of the $j$-th particle, and $\mathcal{L}_t^{(j)}$ the set of all the descendants of $j$ alive at time $t$. For any $x,y$ we have that
\begin{align}\label{625}
	&\P_{\nu}(R_{0,\delta}>x, Q_{0,\delta}<x-y)\nonumber\\
&=\sum_{k=1}^K\P(Z^{(j)}_\delta=0, 1\le j\le k-1, Z^{(k)}_\delta>0, a_k>x, Q_{0,\delta}<x-y)\nonumber\\
	&\le \sum_{k=1}^K \P(Z^{(j)}_\delta=0, 1\le j\le k-1, Z^{(k)}_\delta>0, \inf_{0\le t\le \delta}\max_{u\in\mathcal{L}_t^{(k)}}(\xi^u_t-a_k)<-y)\nonumber\\
	&=\sum_{k=1}^K \P(Z_\delta=0)^{k-1}\P(Z_\delta>0, \inf_{0\le t\le \delta}R_t<-y)\nonumber\\
	&\le \frac{1}{\P(Z_\delta>0)}\P(Z_\delta>0, \inf_{0\le t\le \delta}R_{t}<-y).
	\end{align}
By the many-to-one formula, we have
\begin{align}
\P(Z_\delta>0, \inf_{0\le t\le \delta}R_{t}<-y)&\le \E\left(\sum_{u\in\mathcal{L}_\delta}{\bf 1}_{\{\inf_{0\le t\le \delta}\xi^u_t<-y\}}\ge1\right)\le \E\left(\sum_{u\in\mathcal{L}_\delta}{\bf 1}_{\{\inf_{0\le t\le \delta}\xi^u_t<-y\}}\right)\nonumber\\
&=e^{\lambda \delta}\mP\left(\inf_{0\le t\le \delta}\xi_t<-y\right)\sim e^{\lambda \delta}\frac{q_2}{\alpha}\delta y^{-\alpha}L(y),\qquad y\to\infty,
\end{align}
where the last limit follows from \eqref{inf-xi}.
Combining  the Markov property and  \eqref{625} with $\nu=X_{n\delta}$ and $y=H(n^{-2})$,   we get
		\begin{align}\label{619}
			&\P\left(B_n^c,Q_{n,\delta}<(1-\epsilon)H(e^{-\lambda n\delta}\log n)\left(\vartheta_\delta W_{n\delta}\right)^{\frac{1}{\alpha}}-H(n^{-2}) |\FF_{n\delta } \right)\nonumber\\
			\leq  &\frac{1}{\P(Z_\delta>0)}\P(Z_\delta>0, \inf_{0\le t\le \delta}R_{t}<-H(n^{-2}))
			\le   \frac{ 2\delta e^{\lambda \delta}}{\P(Z_\delta>0)}\frac{q_2}{\alpha}  n^{-2},
		\end{align}
where the last inequality holds for $n$ large enough.
Combining \eqref{419}, \eqref{620} and \eqref{619}, we arrive at \eqref{sum-Q_n}.

Using the Borel-Cantelli lemma and \eqref{sum-Q_n}, we have that for any $\delta>0$,
		$$\liminf\limits_{n\to \infty}\dfrac{Q_{n,\delta} }{H(e^{-\lambda n\delta}\log n)}\ge \left(\vartheta_\delta W\right)^{\frac{1}{\alpha}} \qquad \P^* \textup{-a.s.}$$		
		 Since $H(e^{-\lambda t}\log t)$ is increasing in $t$ for $t$ large, for  any $\delta>0$, we have
		$$\liminf\limits_{t \to \infty}\dfrac{R_t \vee 0}{H(e^{-\lambda t}\log t)} \geq \liminf\limits_{n\to \infty}\dfrac{Q_{n,\delta} \vee 0}{H(e^{-\lambda n\delta}\log n)}\frac{H(e^{-\lambda n\delta}\log n)}{H(e^{-\lambda (n+1)\delta}\log (n+1)\delta )}\ge e^{-\lambda\delta/\alpha}\left(\vartheta_\delta W\right)^{\frac{1}{\alpha}},
$$
$\P^* \textup{-a.s.}$		
		Letting $\delta\to0$, we get the desired result.
		The proof is now complete.
	
		\hfill$\Box$.

It follows from Lemma \ref{cor:M-R} that Lemma \ref{prop:upper-limit} also holds with $M_n$ replaced by $R_n$. To get the limit of $R_t$ as $t\to\infty$, we need to deal with $\sup_{0\le s\le t}R_s$.
If $p_0=0$, then $\sup_{0\le s\le t}R_s=\max_{v\in\mathcal{L}_t}\sup_{0\le s\le t}\xi_s^v.$ Thus by the many-to-one formula, we have that
		\begin{align*}
			\P\left(\sup\limits_{0\leq s \leq t}R_s \ge x\right)\le \E\left( \sum_{v\in\mathcal{L}_t}{\bf 1}_{\{\sup_{0\le s\le t}\xi_s^v\ge x\}}\right) =e^{\lambda t} 	\mP\left(\sup\limits_{0\leq s\leq t}\xi_{s}\ge x\right).
		\end{align*}
		In the following lemma, we show that this assertion is still valid when $p_0>0$.
		
		\begin{lemma}\label{sup-R}
			For any $x>0$ and $t>0$, it holds that
			\begin{align}
				\P\left(\sup\limits_{0\leq s \leq t}R_s \ge x\right) \le   e^{\lambda t}	\mP\left(\sup\limits_{0\leq s\leq t}\xi_{s}\ge  x\right).
			\end{align}		
		\end{lemma}
		
		\noindent {\bf Proof}
		For $r\ge0$ and $x\in\R$, let $\P_{r,x}$ denote the law of $\mathbb{X}$ starting with  an individual  at position $x$ at time $r$, and $\mP_{r,x}$ stands for the law of $\xi$ condition on $\{\xi_r=x\}$. We still use $\sigma_o$ to denote the death time of the initial particle.  Fix $t>0$.
		For any $0\le r<t$ and $x<0$, we have that
		\begin{align}
			\omega(r,x) :=\P_{r,x}\left(\sup\limits_{r\leq s \leq t}R_s< 0\right) &= \P_{r,x}\left(\sup\limits_{r\leq s \leq t}R_s< 0, \sigma_o > t\right) + \P_{r,x}\left(\sup\limits_{r\leq s \leq t}R_s< 0,r \leq \sigma_o \le t\right)  \notag\\
			& =: I_1+I_2 .\notag
		\end{align}
		
		We first consider $I_1$. Let $T_0:=\inf\{u\geq r: \xi_u \ge 0\}$. Since $\sigma_o-r\sim\mathcal{E}(\beta)$ under $\P_{r,x}$,  we have
		\begin{align}
			I_1 &= \mP_{r,x}\left(\sup\limits_{r\leq s \leq t}\xi_s < 0\right)\P_{r,x}(\sigma_o> t) \notag \\
			&= \mP_{r,x}(T_0 > t) e^{-\beta(t-r)} = \mP_{r,x}\left(T_{0} \geq t ,e^{-\int_{r}^{T_{0}\land t}\beta ds}\right) . \notag
		\end{align}
		By the branching property and the strong Markov property, we have
		\begin{align}
			I_2 &= \E_{r,x}\left(\P_{r,x}\left(\sup\limits_{r\leq s \leq t}R_s< 0,r \leq \sigma_o \le t| \mathcal{F}_{\sigma_{o}}\right)\right)    \notag\\
			& = \mE_{r,x}\left(\int_{r}^{t}\mathbf{1}_{\left\{T_{0}> u\right\} } \sum\limits_{k=0}^{\infty} p_k \left(\P_{u,\xi_u}\left(\sup\limits_{u\leq s \leq t}R_s< 0\right)\right)^k\beta e^{-\beta(u-r)} du\right) \notag \\
			& = \mE_{r,x}\left(\int_{r}^{t\land T_{0}}
					f(\omega(u,\xi_u))\beta e^{-\beta(u-r)} du\right), \notag
		\end{align}
		where $f(s)=\sum_{k=0}^\infty p_k s^k$.
		Consequently, $\omega(r,x)$ satisfies the equation
		\begin{align}
			\omega(r,x)= \mP_{r,x}\left(T_{0} > t ,e^{-\int_{r}^{T_{0} \land t}\beta ds}\right)+\mE_{r,x}\left(\int_{r}^{t\land T_{0}} f(\omega(u,\xi_u))\beta e^{-\beta(u-r)} du\right). \notag
		\end{align}
		Then, $1-\omega(r,x)$ satisfies the following equation
		\begin{align}
			1-\omega(r,x)= \mP_{r,x}\left(T_{0} \le  t ,e^{-\int_{r}^{T_{0}\land t}\beta ds}\right)+\mE_{r,x}\left(\int_{r}^{t\land T_{0}}(1- f(\omega(u,\xi_u)))\beta e^{-\beta(u-r)} du\right). \notag
		\end{align}
		Using \cite[Lemma 1.5]{PDE},
		we get that $g(r,x):=1-\omega(r,x)$ also satisfies the following equation
		\begin{align*}
			g(r,x)&= \mP_{r,x}\left(T_{0} \le  t ,e^{\int_{r}^{T_{0}\land t}\lambda ds}\right)+\mE_{r,x}\left(\int_{r}^{t\land T_{0}}(1- f(\omega(u,\xi_u)))\beta e^{\lambda(u-r)} du\right) \\
			& -\mE_{r,x}\left(\int_{r}^{t \land T_0}
			\mu \beta g(u,\xi_u) e^{\lambda(u-r)}du\right), \quad 0\le r<t, x<0 .
		\end{align*}
		\noindent Note that $1-f(1-g)-\mu g \leq 0$. So we have that for $0\le r< t$, and $x<0$,
		\begin{align*}
			g(r,x) \leq  \mP_{r,x}\left(T_{0} \le t ,e^{\int_{r}^{T_{0}\land t}\lambda ds}\right)\le e^{\lambda t}\mP_{r,x}(T_{0} \le t).
		\end{align*}
		Using the time-homogeneity and space-homogeneity of branching L\'evy processes and L\'evy processes, we get that for any $t>0$ and $x>0,$
		\begin{align*}
			\P\left(\sup\limits_{0\leq s \leq t}R_s \ge x\right) & = 1- \P\left(\sup\limits_{0\leq s \leq t}R_s< x\right) =g(0,-x) \notag\\
			& \leq e^{\lambda t}\mP_{0,-x}(T_{0}\le t) =e^{\lambda t} \mP\left(\sup\limits_{0\leq s\leq t}\xi_{s}\ge x\right).
		\end{align*}
		The proof is now complete.
		\hfill$\Box$
		
		\bigskip

		\noindent{\bf Proof of Theorem \ref{them:upper-limit}:}
		Since $\lim_{t\to\infty}\frac{G(t)}{h(t)}=\infty$,
		we have that  $G(t)\overset{t}{\ge}e^{\frac{2\lambda}{3\alpha}t}$. By
		Lemma \ref{cor:M-R} with $3/2<p<2$, we get that
		\begin{align}\label{675}
			\limsup_{t\to\infty}\frac{R_n-M_n}{G(n)}\le0,\quad 	\limsup_{t\to\infty}\frac{M_n-R_n^+}{G(n)}\le0,\quad \P^* \textup{-a.s.}
		\end{align}

		(1) Assume that $\sum_n e^{\lambda n}G(n)^{-\alpha}L(G(n))<\infty$.
		By Theorem \ref{main-theorem-inf}, we have
		$\P^*(\exists T>0, \forall t>T,R_t>0)=1$, which implies that
		$$\limsup_{t\to\infty}\frac{R_t}{G(t)}=\limsup_{t\to\infty}\frac{R_t^+}{G(t)}\ge0, \quad \P^* \textup{-a.s.}.$$
		Let $ V_n:=\sup_{n \le t\le (n+1) }R_t, n\ge0$.
			Since $G$ is non-decreasing, by Lemma \ref{prop:upper-limit} and \eqref{675}, we have, under $\P^*$,
		\begin{align}
			\limsup_{t\to\infty}\frac{R_t}{G(t)}\le &\limsup_{n\to\infty}\frac{V_n}{G(n)}\le \limsup_{n\to\infty}\frac{V_n-R_n}{G(n)}+\limsup_{n\to\infty}\frac{R_n-M_n}{G(n)}+\limsup_{n\to\infty}\frac{M_n}{G(n)}\\
\le &\limsup_{n\to\infty}\frac{V_n-R_n}{G(n)}.
		\end{align}
		By the Borel-Cantelli lemma, to get the desired result, it suffices to prove that, for any $c>0$,
		\begin{align}\label{sum-V}
			\sum_{n=1}^\infty \P^* (V_n-R_n>cG(n ))<\infty.
		\end{align}

		Observe that on $\{Z_n>0\}$, $R_n\in\R$ and
		\begin{align}
			V_n=\max_{u\in\mathcal{L}_n}(\xi^u_n+V_n^u)\le R_n+\max_{u\in\mathcal{L}_n} V_n^u,
		\end{align}
		where $$V_n^u=\sup_{n\le t\le n+1}\max_{v\in\mathcal{L}_t, u\le v}(\xi_t^v-\xi_n^u).$$ It is clear that,  conditioned on $\FF_n$, $\{V_n^u, u\in\cL_n\}$ are i.i.d. with the same law as $(V_0,\P)$.
		Thus
		\begin{align*}
			(1-q)	\P^*(V_n-R_n>cG(n))&\le \P(V_n-R_n>cG(n), Z_n>0)\le \P(\max_{u\in\mathcal{L}_n} V_n^u>cG(n),Z_n>0)\\
			&\le \E(Z_n)\P(V_0>cG(n))=e^{\lambda n}\P(V_0>cG(n))\\
			&\le e^{\lambda n} e^{\lambda}\mP(\sup_{0\le t\le 1}\xi_t>cG(n)),
		\end{align*}
		where the last inequality follows from Lemma \ref{sup-R}. By \eqref{sup-xi},  we have that
		\begin{align}
			e^{\lambda n }\mP\left( \sup_{0\le t\le1}\xi_t\ge cG(n )\right)\sim c^{-\alpha}\frac{q_1}{\alpha} e^{\lambda n}G(n )^{-\alpha}L(G(n )).
		\end{align}
		Thus  \eqref{sum-V} follows immediately.
		
		(2) Assume that $\sum_n e^{\lambda n}G(n)^{-\alpha}L(G(n))=\infty$. By Lemma \ref{prop:upper-limit} and \eqref{675}, we have
		$$\limsup_{t\to\infty}\frac{R_t^+}{G(t)}\ge \limsup_{n\to\infty}\frac{R_n^+}{G(n)}\ge  \limsup_{n\to\infty}\frac{M_n}{G(n)}-\limsup_{n\to\infty}\frac{M_n-R^+_n}{G(n)}=\infty\quad \P^*\mbox{-a.s.},$$
		which implies that $\limsup_{t\to\infty}\frac{R_t}{G(t)}=\infty.$
		
		The proof is now complete.
		$\hfill\qedsymbol$

As a consequence of Theorem \ref{them:upper-limit} and Theorem \ref{main-theorem-inf},  we have the following
result.

\begin{corollary}\label{cor1}
	\begin{itemize}
		\item[(1)]
		\begin{align}\label{1.1}
			\limsup_{t\to\infty}\frac{\log R_t^+-\log H(t^{-1}e^{-\lambda t})}{\log\log t}=\frac{1}{\alpha},   \qquad \P^* \textup{-a.s.}.
		\end{align}
		\item[(2)] $$\lim_{t\to\infty}\frac{\log R_t^+}{t}=\frac{\lambda}{\alpha},\qquad \P^* \textup{-a.s.}$$
	\end{itemize}
	
\end{corollary}

\noindent	{\bf Proof:}
(1) For any $\epsilon>0$,
by Theorem \ref{them:upper-limit} (1) with $G(t)=H(t^{-1}e^{-\lambda t})(\log t)^{(1+\epsilon)/\alpha}$,
 $$\limsup_{t\to\infty}\frac{R_t}{H(t^{-1}e^{-\lambda t})(\log t)^{(1+\epsilon)/\alpha}}=0, \quad \P^* \textup{-a.s.}.$$ which implies that for any  $\delta>0$,
$\frac{R^+_t}{H(t^{-1}e^{-\lambda t})(\log t)^{(1+\epsilon)/\alpha}}\le\delta$ when $t$ is large enough. Thus for $t$ large enough,
$$\log R^+_t-\log H(t^{-1}e^{-\lambda t})\le\log \delta+\frac{1+\epsilon}{\alpha}\log\log t,  \quad \P^* \textup{-a.s.}.$$
Therefore
$$ \limsup_{t\to\infty}\frac{\log R_t^+-\log H(t^{-1}e^{-\lambda t})}{\log\log t}\le \frac{1+\epsilon}{\alpha}, \quad \P^* \textup{-a.s.}.$$
Letting $\epsilon\to 0$, we get $$ \limsup_{t\to\infty}\frac{\log R_t^+-\log H(t^{-1}e^{-\lambda t})}{\log\log t}\le \frac{1}{\alpha}, \quad \P^* \textup{-a.s.}.$$
For any $\epsilon\in(0,1)$,
by Theorem \ref{them:upper-limit} (2) with $G(t)=H(t^{-1}e^{-\lambda t})(\log t)^{(1-\epsilon)/\alpha}$, we have  for any  $\delta>0$,
$\limsup_{t\to\infty}\frac{R^+_t}{H(t^{-1}e^{-\lambda t})(\log t)^{(1-\epsilon)/\alpha}}\ge\delta$, which implies that
$$\limsup_{t\to\infty}\frac{\log R^+_t-\log H(t^{-1}e^{-\lambda t})}{\log\log t}\ge \limsup_{t\to\infty}\frac{\log \delta}{\log\log t}+\frac{1-\epsilon}{\alpha}=\frac{1-\epsilon}{\alpha}.$$
Therefore letting  $\epsilon\to0$,  we get
$$ \limsup_{t\to\infty}\frac{\log R_t^+-\log H(t^{-1}e^{-\lambda t})}{\log\log t}\ge \frac{1}{\alpha},  \quad \P^* \textup{-a.s.}.$$
Hence we have \eqref{1.1}.

(2)
Note that
\begin{align}\label{1.2}
	\log H(t^{-1}e^{-\lambda t})=\frac{\lambda} {\alpha}t+\frac{1}{\alpha}\log t+\log \bar{L}(te^{\lambda t})\sim \frac{\lambda} {\alpha}t.
\end{align}
By \eqref{1.1}, we have
$$\limsup_{t\to\infty}\frac{\log R_t^+}{t}=\frac{\lambda}{\alpha},\qquad \P^* \textup{-a.s.}$$
On the other hand,
by Theorem \ref{main-theorem-inf}, we have
$$\liminf_{t\to\infty} \frac{\log R_t^+}{t}-\frac{\log H(e^{-\lambda t}\log t)}{t}=0,\qquad \P^* \textup{-a.s.}$$
Note that
$$\log H(e^{-\lambda t}\log t)=\frac{\lambda} {\alpha}t-\frac{1}{\alpha}\log \log t+\log \bar{L}(e^{\lambda t}/\log t)\sim \frac{\lambda} {\alpha}t.$$

The proof is now complete.
\hfill$\Box$

\begin{remark}
	Since $H(y)=y^{-1/\alpha}\bar{L}(y^{-1})$, we have
	\begin{align}\label{limit-logH}
		\lim_{y\to0}\frac{	\log H(y)}{\log y^{-1}}=1/\alpha.
	\end{align}
	Assume that $L(x)=(\log x)^r$, where $r\in\R$. By \eqref{limit-logH}, we have that
	$$H(y)=y^{-1/\alpha}(\log H(y))^{r/\alpha}\sim \alpha^{-r/\alpha}y^{-1/\alpha}(\log(y^{-1}))^{r/\alpha},\quad y\to0.$$
	Hence we have
	$$\log H(t^{-1}e^{-\lambda t})-\frac{\lambda} {\alpha}t - \frac{1+r}{\alpha}\log t\to \frac{r}{\alpha}\log(\lambda/\alpha),\quad t\to\infty,$$
	and
	$$ H(e^{-\lambda t}\log t)\sim (\lambda/\alpha)^{r/\alpha}t^{r/\alpha}(\log t)^{-1/\alpha}e^{\lambda t/\alpha},\quad t\to\infty.$$
\end{remark}

	\section{Appendix}	
	In this section, we give some further discussion of
	{\bf Assumption 2}
		We  first recall
		Karamata's theorem (\cite[Theorem 1.5.11]{Bingham}).
	
	\begin{lemma}\label{karamata}
		\begin{itemize}
			\item[(1)]
			If $L$ is slowly varying at $\infty$ and locally bounded in $\left[a,+\infty\right)$ for some $a>0$, then for $r>-1$,
			$$\int_{a}^{x}t^{r}L(t) dt \sim \dfrac{1}{r+1}x^{r+1}L(x) , \quad x\to +\infty.$$
			\item[(2)]
			If $L$ is slowly varying at infinity, then for  $r>1$,
			$$ \int_{x}^{+\infty}t^{-r}L(t)dt \sim \frac{1}{r-1}x^{1-r}L(x ) , \quad x \to +\infty.$$
		\end{itemize}	
	\end{lemma}
	
	\begin{example} {\bf (Strictly Stable process.)} Let $\xi$
		be a strictly $\alpha$-stable process, $\alpha\in (0, 2)$,  on $\R$ with L\'evy measure
		$$n(dy)=c_1x^{-(1+\alpha)}{\bf 1}_{(0,\infty)}(x){\rm d} x+c_2|x|^{-(1+\alpha)}{\bf 1}_{(-\infty,0)}(x){\rm d} x,$$
		where $c_1,c_2\ge 0$, $c_1+c_2\ge0,$
		and if $\alpha=1,$ $c_1=c_2=c$.
		The L\'evy exponent of $\xi$ is given by, for $\theta>0$,
		\begin{align}\label{stable-exponent}
			\psi(\theta)&=\left\{\begin{array}{lll}
				\displaystyle\int_\R (e^{i\theta y}-1-i\theta y) n(\d y), &\alpha\in(1,2);\\
				\displaystyle\int_\R (e^{i\theta y}-1) n(\d y),&\alpha\in (0,1);\\
				\displaystyle \int_\R (e^{i\theta y}-1-i\theta y{\bf 1}_{|y|\le 1}) n(\d y)+i a\theta, &\alpha=1
			\end{array}\right.\nonumber\\
			&=\left\{\begin{array}{lll}
				\displaystyle -\alpha\Gamma(1-\alpha)(c_1e^{-i\pi\alpha/2}+c_2e^{i\pi\alpha/2})\theta^\alpha,&\alpha\in(1,2);\\
				\displaystyle-\alpha\Gamma(1-\alpha)(c_1e^{-i\pi\alpha/2}+c_2e^{i\pi\alpha/2})\theta^\alpha,  &\alpha\in(0,1);\\
				\displaystyle-c\pi\theta+ia\theta, &\alpha=1,
			\end{array}\right.
		\end{align}
		where $a\in\R$ is a constant.
		It is clear that $\psi$ satisfies
		{\bf Assumption 2}.
		For more details on stable processes, we refer
		the readers to \cite[Section 14]{Sato}.

		Note that for $\alpha\neq1$, $e^{-c|\theta|^\alpha}$, $\theta\in \R$,
		is  the characteristic function of a strictly $\alpha$-stable random variable if and only if $|\tan(\pi\alpha/2)|\Re(c)\ge |\Im(c)|$.
	\end{example}
	
	Let $\{(\xi_t)_{t\geq 0}\}$	be a L\'{e}vy process with the generating triplet $(a, b, n)$, that is,
	\begin{align}
		\psi(\theta) &= \log E(e^{i\theta \xi_1}) \notag \\
		&= ia\theta-\frac{1}{2}b^2\theta^2 +\int_{\R/\{0\}}(e^{i\theta y}-1-i\theta y \1_{\{|y|<1\}}) n(dy), \notag
	\end{align}
	where $a\in\R, b\geq 0$, and $\int_{\R/\{0\}}(1\wedge y^2)n(dy) < \infty$.
	Let $r(x):\R \to \R$ be a bounded measurable function, satisfying
	\begin{align}\label{cond:r}
		r(x) = \begin{cases}
			1+o(x),& |x| \to 0; \\
			O(\frac{1}{|x|}) , & |x| \to \infty.
		\end{cases}
	\end{align}
	Then $\psi$ can be rewritten as
	\begin{align*}
		\psi(\theta)= ia_r\theta-\frac{1}{2}b^2\theta^2 +\int_{\R/\{0\}}(e^{i\theta y}-1-i\theta y r(y)) n(dy) ,
	\end{align*}
	where $a_r=a+\int_{\R/\{0\}}y(r(y)-\1_{\{|y|<1\}})n(dy)$. The triplet of $\xi$ is written as $(a_r,b,n)_r$.

For any $t>0$, we define
a measure $n_t(dy)$ as follows: for any positive function $g$,
$$ \int_{\R/\{0\}}g(y)n_t(dy) =\int_{\R/\{0\}}g(yt^{-1})
n(dy).$$

		\begin{lemma}\label{lem1}
		If there exist $q_1, q_2\ge0$ such that
		\begin{align}\label{lim-nt}
			\lim\limits_{t \to +\infty}t^{\alpha}L(t)^{-1}n(t,+\infty) =\dfrac{q_1}{\alpha}, \mbox{ and } \lim\limits_{t \to +\infty}t^{\alpha}L(t)^{-1}n(-\infty,-t) =\dfrac{q_2}{\alpha},
		\end{align}
		then as $\theta \to 0^+$,
		$$\theta^{-\alpha}L(\theta^{-1}) \int_{\R} (e^{i\theta y}-1 -i\theta y \1_{\{|\theta y|<1\}})n(dy)\to  \int_{\R} (e^{iy}-1-iy\1_{\{|y|<1\}}) \nu_\alpha(dy), $$
		where $\nu_{\alpha}(dx) := q_1x^{-1-\alpha}\1_{(0,+\infty)}(x)dx+q_2|x|^{-1-\alpha}\1_{(-\infty,0)}(x)dx.$
		
		Furthermore, if $0<\alpha<1$, then as $\theta \to 0^+$,
		\begin{align}\label{7.4}
		\theta^{-\alpha}L(\theta^{-1})\int_{\R} (e^{i\theta y}-1-i\theta y{\bf 1}_{\{|y|<1\}})n(dy) \to  -\alpha\Gamma(1-\alpha)(q_1e^{-i\pi\alpha/2}+q_2e^{i\pi\alpha/2});
		\end{align}
		if  $1<\alpha<2$, then as $\theta \to 0^+$,
		\begin{align}\label{7.5}
		\theta^{-\alpha}L(\theta^{-1})\int_{\R} (e^{i\theta y}-1-i\theta y)n(dy) \to  -\alpha\Gamma(1-\alpha)(q_1e^{-i\pi\alpha/2}+q_2e^{i\pi\alpha/2}).
		\end{align}
		
	\end{lemma}
	\noindent{\bf{Proof:}} 	
For any $t>0$, let $\tilde{n}_t(dy):=t^\alpha L(t)^{-1}n_t(dy)$. By \eqref{lim-nt}, we have that for any $x>0$,$$\tilde{n}_t(x,\infty)=t^\alpha L(t)^{-1}n(tx,\infty)\to \frac{q_1}{\alpha}x^{-\alpha}=\nu_\alpha(x,\infty)$$
and
$$\tilde{n}_t(-\infty,-x)=t^\alpha L(t)^{-1}n(-\infty,-tx)\to \frac{q_2}{\alpha}x^{-\alpha}=\nu_\alpha(-\infty,-x).$$
Thus, for any $g\in C_b^0(\R)$,
\begin{align}\label{7.7}
\int_0^\infty g(y)\tilde{n}_t(dy)\to \int_0^\infty g(y)\nu_\alpha(dy).
\end{align}

We claim that
	\begin{align}\label{7.1}
	\lim\limits_{\varepsilon \downarrow 0}\limsup\limits_{t\to +\infty} \int_{0<|y|<\varepsilon}y^{2} \widetilde{n_t}(dy) =0.
	\end{align}
	In fact, by Fubini's theorem, we have that
	\begin{align*}
		\int_{0<y<\varepsilon}y^2 \widetilde{n_t}(dy) &= t^{\alpha}L(t)^{-1} \int_0^\epsilon y^2 n_t(dy) = t^{\alpha}L(t)^{-1} \int_0^\epsilon\left(\int_{0}^{y}2x dx\right) n_t(dy) \\
		& = t^{\alpha}L(t)^{-1} \int_{0}^{\varepsilon}2x dx \int_{x}^{\varepsilon}n_t(dy)\le t^{\alpha}L(t)^{-1} \int_{0}^{\varepsilon}2xn_t(x,\infty)dx \\
		&=t^{\alpha-2}L(t)^{-1} \int_{0}^{t\varepsilon}2xn(x,\infty)dx.
	\end{align*}	
	Note that $$t^{\alpha-2}L(t)^{-1} \int_{0}^{1}2xn(x,\infty)dx=t^{\alpha-2}L(t)^{-1}\int(1\wedge y^2)n(dy)\to0,\quad t\to\infty.$$
	By \eqref{lim-nt} and
	 Karamata's theorem, we have, as $t\to\infty$,
	$$ \int_{1}^{t\epsilon}xn(x,\infty)dx\sim\frac{q_1}{\alpha}\frac{(t\epsilon)^{2-\alpha}L(t)}{2-\alpha}.$$
	Thus we have
	$$ \limsup\limits_{t\to +\infty} \int_{0}^{\varepsilon}y^{2} \widetilde{n_t}(dy) \le \frac{q_1}{\alpha}\dfrac{\varepsilon^{2-\alpha}}{2-\alpha}.$$
	Similarly, we have
	$$ \limsup\limits_{t\to +\infty} \int_{-\varepsilon}^0 y^{2} \widetilde{n_t}(dy) \le \frac{q_2}{\alpha}\dfrac{\varepsilon^{2-\alpha}}{2-\alpha}.$$
Letting $\epsilon\to0$, we get \eqref{7.1}.

	By \eqref{7.7} and \eqref{7.1}, and applying  \cite[Theroem 8.7]{Sato}, we obtain that, for any bounded continuous function $r$ satisfying \eqref{cond:r},
	\begin{align}\label{58}
		\lim\limits_{t \to +\infty}\int_{\R} (e^{i\theta y}-1-i\theta y r(y))\widetilde{n_t}(dy)= \int_{\R}(e^{i\theta y}-1-i\theta yr(y))\nu_\alpha(dy).
	\end{align}
Since $v_\alpha(\{y\})=0$ for all $y\in\R$, thus  \eqref{58} holds for $r(y)={\bf 1}_{\{|y|\le 1\}}$.
Thus we have that
	\begin{align}\label{510}
		\lim\limits_{t \to +\infty}t^{\alpha}L(t)^{-1}\int_{\R} (e^{it^{-1} y}-1-it^{-1} y \1_{\{|y|< t\}})n(dy)= \int_{\R}(e^{i y}-1-i\ y\1_{\{|y|< 1\}})\nu_{\alpha}(dy).
	\end{align}

(1) Now we assume that  $0<\alpha<1$.
 Note that for $t>1$,
	\begin{align*}
		\int_{1}^{t}yn(dy)&=\int_{1}^{t}\left(\int_{0}^{y}dx \right)n(dy)=\int_{0}^{1}dx\int_{1}^{t}n(dy)+\int_{1}^{t}dx\int_{x}^{t}n(dy)\\
		&=n(1,t)+\int_{1}^{t}n(x,t)dx\\ &=n(1,t)+\int_{1}^{t}n(x,+\infty)dx-(t-1)n(t,+\infty).
	\end{align*}
	By Lemma \ref{karamata}, we have that $\int_{1}^{t}n(x,+\infty)dx\sim \frac{q_1}{\alpha(1-\alpha)}t^{1-\alpha}L(t)$. Thus we have
	\begin{align}\label{7.3}
		\lim\limits_{t \to +\infty}t^{\alpha}L(t)^{-1}t^{-1}\int_{1}^{t}yn(dy)=\dfrac{1}{1-\alpha}\dfrac{q_1}{\alpha}-\dfrac{q_1}{\alpha}= \int_{0}^{1}y\nu_\alpha (dy).
	\end{align}
	Similarly, we have
	\begin{align}\label{7.6}
		\lim\limits_{t \to +\infty}t^{\alpha}L(t)^{-1}t^{-1}\int_{-t}^{-1}yn(dy)=\dfrac{1}{1-\alpha}\dfrac{q_2}{\alpha}-\dfrac{q_2}{\alpha}= \int_{-1}^{0}y\nu_\alpha (dy).
	\end{align}
	By \eqref{510}, \eqref{7.3} and \eqref{7.6}, we obtain that \eqref{7.4} holds.
	
(2) Now we consider the case when $\alpha\in(1,2)$.	
	By Fubini's theorem and lemma \ref{karamata}, we have that
	\begin{align}
		\int_{t}^{+\infty}y n(dy)&=\int_{t}^{+\infty}\left( \int_{0}^{y}dx\right)n(dy)=\int_0^t dx\int_t^\infty n(dy)+\int_{t}^{+\infty}dx\int_{x}^{+\infty}n(dy)\notag \\
		&=tn(t,+\infty)+\int_{t}^{+\infty}n(x,+\infty)dx\notag \\
		& \sim \frac{q_1}{\alpha}t^{\alpha-1}L(t)+\frac{q_1}{\alpha(\alpha-1)}t^{\alpha-1}L(t)=\frac{q_1}{\alpha-1}t^{\alpha-1}L(t).
	\end{align}
	Thus
	\begin{align}\label{511}
	 \lim\limits_{t \to +\infty}t^{\alpha}L(t)^{-1}\int_{t}^{+\infty} t^{-1}yn(dy)= \int_{1}^{+\infty}y\nu_\alpha(dy).
	\end{align}
	Similarly, we have
	\begin{align}\label{511'}
		\lim\limits_{t \to +\infty}t^{\alpha}L(t)^{-1}\int_{-\infty}^{-t} t^{-1}yn(dy)=\int_{-\infty}^{-1}y\nu_\alpha(dy).
	\end{align}
Thus by \eqref{510}, \eqref{511} and \eqref{511'}, we obtain that \eqref{7.5} holds.

The proof is now complete.
	$\hfill\qedsymbol$	
	
	\begin{theorem}
	There exist a constant $c_*$, $\alpha \in (0,2)$ and a function $L$  slowly varying at $\infty$ such that
		\begin{align}\label{51}
			\psi(\theta) \sim -c_*\theta^{\alpha}L(\theta^{-1}) , \theta \to 0^+,
		\end{align}
	 if and only if the following two conditions hold:
		\begin{itemize}
			\item[(1)] there exists nonnegative numbers $q_1,q_2$ such that
	\begin{align*}
			\lim\limits_{t \to +\infty}t^{\alpha}L(t)^{-1}n(t,+\infty) =\dfrac{q_1}{\alpha}, \quad \mbox{and}\lim\limits_{t\to +\infty}t^{\alpha}L(t)^{-1}n(-\infty,-t) =\dfrac{q_2}{\alpha};
		\end{align*}

		\item[(2)] if $\alpha\in(1,2)$,  $a+\int_{\{|y|>1\}}yn(dy)=0$;
		if $\alpha=1$, $q_1=q_2$ and   $$\lim_{t\to\infty}L(t)^{-1}\left(a+\int_{1<|y|<t} yn(dy) \right)
		=:c_0
		$$
		exists.
			\end{itemize}
			Moreover, the relationship between $c_*$ and  $(q_1,q_2,c_0)$ is as follows:
			if $\alpha\neq 1$,
			$$c_*=\alpha\Gamma(1-\alpha)\left(q_1e^{-i\pi\alpha/2}+q_2e^{i\pi\alpha/2}\right),$$
			and if $\alpha=1$, $c_0=-\Im(c_*)$ and
			$q_1=q_2=\Re(c_*)/\pi.$
	\end{theorem}
\noindent	{\bf Proof:}
First, we prove the sufficiency. Assume the two conditions hold. If  $\alpha\in(0,1)$, then by  Lemma \ref{lem1}, we have that as $\theta\to 0^+$,
$$\theta^{-\alpha}L(\theta^{-1})^{-1}\psi(\theta)\to -\alpha\Gamma(1-\alpha)\left(q_1e^{-i\pi\alpha/2}+q_2e^{i\pi\alpha/2}\right).$$
If $\alpha\in(1,2)$, then as $\theta\to 0^+$,
$$\psi(\theta)=-\frac{1}{2}b^2\theta^2+\int_{\R} (e^{i\theta y}-1-i\theta y)n(dy)\sim -\alpha\Gamma(1-\alpha)\left(q_1e^{-i\pi\alpha/2}+q_2e^{i\pi\alpha/2}\right)\theta^{\alpha}L(\theta^{-1}).$$
If $\alpha=1$, then as  $\theta\to 0^+$,
\begin{align*}
\psi(\theta)&=i\left(a+\int_{1<|y|<\theta^{-1}} yn(dy) \right)\theta--\frac{1}{2}b^2\theta^2+\int_{\R} (e^{i\theta y}-1-i\theta y{\bf 1}_{\{|y|\le \theta^{-1}\}})n(dy)\\
&\sim (ic_0-q_1\pi)\theta L(\theta^{-1}).
\end{align*}

Now we assume that \eqref{51} holds.
 It is clear that \eqref{51} is equivalent
	\begin{align}\label{54}
		\lim\limits_{t \to +\infty}e^{t^{\alpha}L(t)^{-1}\psi(\theta t^{-1})} =
		e^{\widetilde{\psi}(\theta)}, \quad \theta>0,
	\end{align}
	where
	$\widetilde{\psi}(\theta)=-c_*\theta^\alpha.$

	Note that the left side of \eqref{54} is  the  characteristic
	function
	of an infinitely divisible random variable $Y_t$ with L\'evy measure $\tilde{n}_t(dy)$.
By \cite[Theorem 8.7 (1)]{Sato}, if \eqref{54} holds, then
	$e^{\tilde{\psi}(\theta)}$ is the characteristic
	function of an infinitely divisible random variable. By the expression of $\tilde{\psi}(\theta)$, $e^{\widetilde{\psi}(\theta)}$ must be
	the characteristic function of a strictly $\alpha$-stable random variable $Y$. Thus if $\alpha\neq1$, then $|\tan(\pi\alpha/2)|\Re(c_*)\ge |\Im(c_*)|$. Consequntely, the L\'evy measure of $Y$ is given by
	$$\nu_{\alpha}(dx) := q_1x^{-1-\alpha}\1_{(0,+\infty)}(x)dx+q_2|x|^{-1-\alpha}\1_{(-\infty,0)}(x)dx,$$
where $q_1\ge 0$ and $q_2\ge0$  satisfy the  following equation:
if $\alpha\neq 1$,
$$c_*=\alpha\Gamma(1-\alpha)\left(q_1e^{-i\pi\alpha/2}+q_2e^{i\pi\alpha/2}\right),$$
and if $\alpha=1$,
$q_1=q_2=\Re(c_*)/\pi.$

By \cite[Theorem 8.7 (1)]{Sato}, we get that  for any $g \in C_b^0(\R)$,
		\begin{align}\label{52}
		\lim\limits_{t \to +\infty}t^{\alpha}L(t)^{-1} \int_{\R/\{0\}} g(t^{-1}y)n(dy) = \int_{\R/\{0\}} g(y) \nu_{\alpha}(dy).
	\end{align}
	Because   $\nu_{\alpha}(\{x\}) =0 ,\forall x \in \R $   ,then \eqref{52} holds for $g(x)={\bf 1}_{(1,\infty)}(x)$ and $g(x)={\bf 1}_{(-\infty,-1)}(x)$, i.e.,
	\begin{align*}
		\lim\limits_{t \to +\infty}t^{\alpha}L(t)^{-1}n(t,+\infty) = \nu_{\alpha}(1,+\infty)=\dfrac{q_1}{\alpha};\\
		\lim\limits_{t\to +\infty}t^{\alpha}L(t)^{-1}n(-\infty,-t) = \nu_{\alpha}(-\infty,-1)=\dfrac{q_2}{\alpha}.
	\end{align*}
Now using Lemma \ref{lem1}, we can get the second condition holds.

	The proof is now complete.
	$\hfill\qedsymbol$	
	\bigskip

	\begin{remark}
		For $\alpha =1 $, we can assume that L\'evy measure $n$ is symmetric and $a=0$.
	\end{remark}
	
	\section*{Acknowledgments}
	Yan-Xia Ren  is supported by   NSFC (Grant No. 12231002) and the Fundamental Research Funds for Central Universities, Peking University LMEQF.
	Renming Song's research  was supported in part by a grant from the Simons
	Foundation (\#960480, Renming Song). Rui Zhang is  supported by  NSFC (Grant No.  12271374, 12371143).

\end{document}